\documentclass[12pt]{amsart}

\usepackage{amsfonts,amsmath,latexsym,amssymb,verbatim,amsbsy,amsthm}

\usepackage{graphicx}
\usepackage{nicefrac}
\usepackage{dsfont}
\usepackage{mathrsfs}
\usepackage{cancel}
\usepackage[normalem]{ulem}
\usepackage{color}
\usepackage{enumitem}

\usepackage[top=1in, bottom=1in, left=1in, right=1in]{geometry}

\usepackage[dvipsnames]{xcolor}
\usepackage{relsize}

\usepackage[colorlinks=true, pdfstartview=FitV, linkcolor=RoyalBlue,citecolor=red, urlcolor=blue]{hyperref}

\numberwithin{equation}{section}
\numberwithin{figure}{section}

\renewcommand{\geq}{\geqslant}

\renewcommand{\leq}{\leqslant}

\newcommand{\ds}{\displaystyle} 
\newcommand{\be}{\begin{equation}}
\newcommand{\ee}{\end{equation}}

\theoremstyle{plain}
\newtheorem{THEOREM}{Theorem}[section]

\newtheorem{theorem}[THEOREM]{Theorem}
\newtheorem{corollary}[THEOREM]{Corollary}
\newtheorem{lemma}[THEOREM]{Lemma}
\newtheorem{proposition}[THEOREM]{Proposition}

\theoremstyle{plain} 

\newtheorem{definition}[THEOREM]{Definition}

\theoremstyle{remark}

\newtheorem{remark}[THEOREM]{Remark}

\newcommand{\myr}[1]{{{#1}}} 
\newcommand{\myout}[1]{}
\newcommand{\mycancel}[1]{}
\newcommand{\myangle}[1]{(1+{#1})} 

\newcommand{\rem}[1]{Remark~\ref{#1}}

\newcommand\widebar[1]{{\mathop{\overline{#1}}}}

\newcommand\argmax{{\mathop{\textnormal{argmax}}}}

\def\bu{{\mathbf u}}
\def\bw{{\mathbf w}}

\def\bm{{\mathbf m}}

\def\rhoa{\rho} 
\def\rhob{\rho} 
\def\rhobp{\rho'} 
\def\ea{e_{\!{}_\Pressure}} 
\def\eb{e_{\!{}_\Pressure}} 
\def\ebp{e'_{\!{}_\Pressure}} 
\def\bua{\bu} 
\def\bub{\bu} 
\def\bubp{\bu'} 
\def\Pressurea{\Pressure} 
\def\bqa{\bq} 
\def\aligna{\align} 
\def\fa{f_N} 
\def\fb{f_N} 
\def\fbp{f'_N} 
\def\Ea{E} 
\def\Ma{M} 
\def\bma{\bm} 
\def\phiab{\phi} 
\def\kab{k} 
\def\Na{N} 
\def\Nb{N} 
\def\na{{\mathbf n}} 
\def\Qa{Q} 
\def\Qp{Q_p}
\def\Da{D} 
\def\Omegaa{{\mathcal S}} 
\def\Omegab{{\mathcal S}} 
\def\Omegae{{}}

\def\bubar{\widebar{\bu}}
\def\bvbar{\widebar{\bv}}
\def\ubar{\widebar{u}}

\def \bx {{\mathbf x}}

\def \bxp {{\mathbf x}'}
\def \by {\bxp} 
\def \xp {{x}'}
\def \bu {{\mathbf u}}

\def \bv {{\mathbf v}}
\def \bvp {{\mathbf v}'}
\def \bz {{\mathbf z}}

\def\Etotal{{\mathscr E}} 

\def\delE{\delta\Etotal}
\def\delbu{\delta \bu} 
\def\delbv{\delta \bv} 

\newcommand{\R}{\ensuremath{\mathbb{R}}}   

\def \lan {\langle}
\def \ran {\rangle}

\def\d{\textnormal{d}}    

\def \hf {\frac{1}{2}}

\def \dx  {\, \d\bx} 
\def\dxp {\, \d\bxp} 
\def\dvp {\, \d\bv'} 
\def \dy  {\dxp} 
\def \dz  {\, \mathrm{d}\bz}

\def \dv  {\, \mathrm{d}\bv}

\def \ddt  {\frac{\d}{\d t}} 

\def \Pressure {{\mathbb P}}
\newcommand{\pressure}{\hspace*{0.1ex}{\raisebox{.4ex}{${}_{\mathbb P}$}}\hspace*{-0.1ex}} %
\def \align {{\mathbf A}}

\newcommand{\bq}{\textbf{q}}

\begin{document}

\title[Swarming: hydrodynamic alignment with pressure]{Swarming: hydrodynamic alignment with pressure}

\author{Eitan Tadmor}
\address{Department of Mathematics and Institute for Physical Sciences \& Technology (IPST), University of Maryland, College Park}
\email{tadmor@umd.edu}

\date{\today}



\thanks{\textbf{Acknowledgment.} Research was supported by ONR grant N00014-2112773.}

\begin{abstract}
We study the swarming behavior of hydrodynamic alignment. Alignment reflects  steering towards a weighted average heading. We consider  the  class of so-called $p$-alignment hydrodynamics, based on $2p$-Laplacians, and weighted by  a general family of  symmetric communication kernels.
 The main new aspect here is the long time emergence behavior for a general class of pressure tensors \emph{without} a closure assumption, beyond the mere requirement  that they form an energy dissipative process. We refer to such pressure laws as `entropic', and    prove the flocking  of $p$-alignment hydrodynamics, driven by singular kernels with general class of entropic pressure tensors. These  results indicate the rigidity of alignment in driving long-time flocking behavior  despite the lack of thermodynamic  closure.
\end{abstract}

\maketitle
\setcounter{tocdepth}{1}
\tableofcontents

\section{Introduction --- alignment dynamics and entropic pressure}\label{sec:entropic}
Alignment reflects  steering towards average heading, \cite{Rey1987}. It plays an indispensable role in the process of emergence  in swarming dynamics, and in particular --- in flocking, herding,  schooling,..., \cite{VCBCS1995,CF2003,CKFL2005,CS2007a,CS2007b,Bal2008,Kar2008,VZ2012,MCEB2015,PT2017}, as well as the formation of other self-organized clustering in human interactions and in dynamics of sensor-based networks,  \cite{Kra2000,BeN2005,BHT2009,JJ2015,RDW2018,DTW2019,Alb2019}; more can be found in \cite[\S9]{MT2014}, in the book series on active matter,  \cite{BDT2017/19, BCT2022}, and in the recent Gibbs' lecture \cite{Tad2022a}.\newline 
We discuss alignment dynamics in two parallel descriptions. Historically, alignment models were introduced in the context of agent-based description \cite{Aok1982,Rey1987,VCBCS1995}. In particular, our discussion is motivated by the celebrated Cucker-Smale model, \cite{CS2007a,CS2007b}, in which alignment is governed by  weighted graph Laplacians. Our main  focus, however, is on the corresponding  hydrodynamic description, the so-called  Euler alignment  equations, governed by  a general class of weighted $p$-graph Laplacians, \cite{HT2008, CFTV2010, HHK2010, Shv2021}. In both cases --- the agent-based and hydrodynamic descriptions, the weights for the protocol of alignment reflect  pairwise interactions,  and are quantified by proper \emph{communication kernel}.  Communication  kernels are either derived empirically, deduced from higher-order principles, learned from the data, or
postulated based on phenomenological arguments, e.g.,
\cite{CS2007a,CDMBC2007, Bal2008,GWBL2012, JJ2015,LZTM2019,MLK2019, ST2020b}. The specific structure of such kernels, however, is not necessarily known. Instead, we ask how different classes of communication kernels affect the swarming behavior.\newline
 The passage from agent-based to hydrodynamic descriptions requires a proper notion of hydrodynamic pressure. In section \ref{sec:entropic}  we introduce a  class of \emph{entropic pressures} for hydrodynamic alignment and in section \ref{sec:p-alignment} we extend the discussion to the larger class of hydrodynamic $p$-alignment. Our goal is to make a systematic study of the long-time swarming behavior of hydrodynamic alignment, portrayed in section \ref{sec:swarming}, with  entropy pressure laws. Specifically, we use the decay of \emph{energy fluctuations}, discussed in section \ref{sec:fluctuations}, in order to quantify the emergence of flocking behavior, depending on the communication kernel. Almost all available literature is devoted to the case of `pressure-less' alignment. We review these results   in section \ref{sec:mono-kinetic}. The main theme here is unconditional flocking for pressure-less $p$-alignment, driven by \emph{heavy-tailed} communication kernels. In section \ref{sec:with-pressure} we discuss  hydrodynamic alignment  driven by a general class of entropic pressure. The remarkable aspect here is that despite the lack of closure of such entropic  pressure laws, there holds unconditional flocking of $p$-alignment driven  by \emph{singular, heavy-tailed}  communication kernels.  We are aware that the methodology developed here can be utilized with other Eulerian-based dissipative systems.\newline
 The detailed computations are outlined in appendix \ref{sec:derivation},  \ref{sec:pointwise}, \ref{sec:GN} and \ref{sec:dispersion}.

\subsection{Hydrodynamic description of alignment}
We study the long-time behavior of the    (hydro-)dynamic  description for alignment,
\begin{subequations}\label{eqs:hydro}
\begin{equation}\label{eq:hydro}
    \left\{\begin{array}{c}
    \begin{split}
    & \partial_{t}\rhoa+\nabla_\bx\cdot(\rhoa\bua) = 0,\\
    & \partial_{t}(\rhoa\bua)+\nabla_\bx\cdot(\rhoa\bua\otimes\bua+\Pressurea) = \aligna(\rho,\bu),
    \end{split}      
    \end{array}\right. \quad (t,\bx)\in (\R_t, \R^d).
\end{equation}
The dynamics is captured by density $\rhoa: \R_t\times \R^d \mapsto \R_+$, momentum, $\rho\bu: \R_t\times \R^d \mapsto \R^d$, and pressure tensor, $\Pressurea: \R_t\times \R^d \mapsto \R^d\times \R^d$,
subject to initial data $\displaystyle (\rhoa,\bua,\Pressurea)_{|_{t=0}} = (\rho_{0},\bu_{0},\Pressure_{0})$, and is driven by an
\emph{alignment} term acting on the support  ${\mathcal S}(t):=\textnormal{supp}\,\rho(t,\cdot)$,
\begin{equation}\label{eq:align}
        \aligna(\rho,\bu) := \int \limits_{{\mathcal S}(t)}\phiab (\bx,\by)(\bub(t,\by)-\bua(t,\bx))\rhoa(t,\bx)\rhob(t,\by)\dy, \qquad \phiab (\bx,\by)=\phiab (\by,\bx).
\end{equation} 
The alignment term on the right reflects  steering towards average heading.
 Here, different weighted averages are dictated by  symmetric communication kernels $\phiab(\cdot,\cdot)$. Prototypical examples include  \emph{metric kernels}, 
$\phiab (\bx,\by)=\kab(|\bx-\by|)$,
which go back to \cite{CS2007a}. Other classes of symmetric kernels that either dictated by the problem or learned from the data can be found in \cite{GWBL2012,JJ2015,LZTM2019,MLK2019}, and finally we mention 
 topologically-based kernels studied in \cite{ST2020b}, 
$\phiab (\bx,\by)=\kab(m(C(\bx,\by)))$,
 where $\ds m(C(\bx,\by))= \int_{C}\rho(t,\bz)\dz$ is the  mass enclosed in an intermediate domain $C=C(\bx,\bxp)$ with tips at $\bx$ and $\by$. The prominent role of metric kernels enters when we assume that there exists a  radial kernel, $\kab(r)$, such that 
\begin{equation}\label{eq:radial}
\phiab (\bx,\by)\geq \kab (|\bx-\by|).
\end{equation}
\end{subequations}
We further assume that the metric kernel $\kab (r)$ is decreasing with the distance $r$, reflecting the typical observation that the intensity of alignment decreases with the distance. In particular, we address general  metric kernels $\phiab (|\cdot|)$ whether decreasing or not,  in terms of their \emph{decreasing envelope} $\kab (r):=\min\{\phiab (|\bx|) \ | \ |\bx|\leq r\}$.  Observe that we do not place any restriction on the upper-bound of $\phi$; in particular, therefore, our discussion includes the important sub-class of \emph{singular} communication kernels $\kab(r)=r^{-\alpha}, \ \alpha>0$, \cite{ST2017a, DKRT2018,MMPZ2019,AC2021b}.
\subsection{Entropic pressure}
System \eqref{eqs:hydro}  is not closed in the sense that the pressure $\Pressure$ is not specified ---  neither in terms of algebraic   relations with $(\rhoa,\bu)$, nor do we specify  a precise dynamics of $\Pressure$. Indeed,  We do not dwell here on the details of the  underlying the pressure tensor.  Instead, we treat a rather general class  of  pressure laws satisfying an essential structural (dissipative) property which, as we shall show, maintains  long time flocking behavior. This brings us to the following.
\begin{definition}[{\bf Entropic pressure}]\label{def:meso-pressure}
We say that $\Pressurea$ is \myr{an} entropic pressure  associated with \eqref{eqs:hydro} if it has a non-negative  trace, $\rhoa \ea:=\hf \text{trace}(\Pressurea)\geq0$, which satisfies
\begin{equation}\label{eq:meso-pressure}
    \partial_{t}(\rhoa \ea )+\nabla_{\bx}\cdot(\rhoa \ea \bua +\bqa )+\textnormal{trace}(\Pressurea\nabla\bua) \leq -   2\int \limits_{{\mathcal S}(t)} \phiab (\bx,\bxp)\ea(t,\bx)\rho(t,\bx)\rho(t,\bxp) \dxp.
\end{equation}
Here $\bqa$ is an arbitrary $C^1$-flux. 
\end{definition}

\noindent
{\bf Why entropic pressure?} System \eqref{eqs:hydro} falls under the general category of hyperbolic balance laws \cite[Chapter III]{Daf2016}, and \eqref{eq:meso-pressure} can be viewed as  an entropy inequality associated with such balance law.  To this end, we note that a formal manipulation  of the mass and momentum equations,  \eqref{eq:hydro}$\ds {}_1 \times \frac{|\bu|^2}{2} + $  \eqref{eq:hydro}${}_2 \cdot \bu$ yields\footnote{Here and below for a quantity $\square=\square(t,\bx)$ we abbreviate $\square':=\square(t,\by)$} 
 \begin{equation}\label{eq:energy-kin}
  \partial_t \Big(\frac{\rhoa}{2}|\bua|^2\Big) + \nabla_\bx\cdot \Big(\frac{ \rhoa}{2}|\bua|^2\bua +\Pressurea\bua\Big) - \textnormal{trace} \big(\Pressurea \nabla\bua\big)  = 
 -  \int  \limits_{{\mathcal S}(t)}\phiab (\bx,\by)\rhoa\bua\cdot(\bua-\bubp)\rhobp \dy.
 \end{equation}
Adding  the entropic description of the pressure postulated in \eqref{eq:meso-pressure} leads to the entropic statement for  the total energy, $\ds \Ea:= \frac{|\bua|^2}{2} +\ea$, 
\begin{equation}\label{eq:entropy-ineq}
\partial_t (\rho E) + \nabla_\bx\cdot(\rho E\bu+\Pressure\bu+\bq) \leq 
 -\int \limits_{{\mathcal S}(t)}\phiab (\bx,\by)\big(|\bu|^2-\bua\cdot\bubp+2\ea\big)\rho\rhobp \dy. 
\end{equation}
Thus, the notion of entropic pressure \eqref{eq:meso-pressure} complements the 
balance laws in \eqref{eqs:hydro} to form the \emph{entropy inequality} \eqref{eq:entropy-ineq}.

To further motivate why this notion of an entropic pressure, we appeal to its underlying \emph{kinetic formulation}. The hydrodynamics \eqref{eqs:hydro} corresponds to  the large-crowd dynamics 
  of $\Na$ agents with position/velocity  $\displaystyle (\bx_{i}(t),\bv_{i}(t)): \R_t\mapsto \R^{d}\times\R^{d}$, governed by the celebrated agent-based  alignment model of Cucker \& Smale \cite{CS2007a,CS2007b}
\begin{equation}\label{eq:CS}
    \left\{\begin{array}{c}
    \begin{split}
        \ddt\bx_{i}(t) &= \bv_{i}(t),   \\
        \ddt \bv_{i}(t) &= \frac{1}{\Nb}\sum_{j = 1}^{\Nb}\phiab_{ij}(t)(\bv_{j}(t)-\bv_{i}(t)),
    \end{split}      
    \end{array}\right. \qquad i=1,2,\ldots \Na.
\end{equation}
The alignment dynamics is driven by a weighted graph Laplacian on the right of \eqref{eq:CS}${}_2$, dictated by the symmetric communication kernel,
$\phi_{ij}(t):=\phiab(\bx_{i}(t),\bx_{i}(t))$.
 The passage from the agent-based to the hydrodynamic description is realized by moments of the \emph{empirical distribution} 
\begin{equation*}
    \fa(t,\bx,\bv) := \frac{1}{\Na  }\sum_{i=1}^{\Na  }\delta_{\bx_{i}(t)}\otimes\delta_{\bv_{i}(t)}, \qquad (t,\bx,\bv)\in \R_{t}\times\R^d\times\R^{d}.
\end{equation*}
The large crowd limits which are assumed to exist,  recover \eqref{eqs:hydro} 
with
\[
\rhoa(t,\bx)= \lim_{\Na  \rightarrow\infty}\int \limits_{\R^d} \fa(t,\bx,\bv)\dv \ \ \textnormal{and} \ \  
\rhoa\bua(t,\bx) = \lim_{\Na  \rightarrow\infty}\int \limits_{\R^d} \bv \fa(t,\bx,\bv)\dv.
\] 
This passage from agent-based to macroscopic description is outlined  in appendix \ref{sec:hydro-des} below. It was justified for smooth kernels \cite{HT2008,CFTV2010,CCR2011,FK2019,NP2021,Shv2021} and   at least  mildly singular kernels, \cite{Pes2015,PS2019,MMPZ2019}.  In this context, the pressure or Reynolds stress tensor   corresponds to the \emph{second-order moments}
\begin{equation}\label{eq:pressure}
        \Pressurea(t,\bx) = \lim_{\Na  \rightarrow\infty}\int \limits_{\R^d} (\bv-\bua)\otimes (\bv-\bua)\fa(t,\bx,\bv)\dv.
\end{equation}
  We observe that the kinetic description of pressure in \eqref{eq:pressure} is consistent with the entropic inequality postulated in \eqref{eq:meso-pressure}.    Indeed,  $\rhoa\ea:=\hf \text{trace}(\Pressurea)$ is the \emph{internal energy} which quantifies  microscopic  fluctuations around the bulk velocity $\bua$,
\begin{equation}\label{eq:internale}
\rhoa \ea = \lim_{\Na  \rightarrow\infty} \int \limits_{\R^d} \hf |\bv-\bua|^2\fa(t,\bx,\bv)\dv.
\end{equation}
This kinetic description of internal energy  yields (detailed derivation is carried   out in appendix \ref{sec:meso-ppressure} below), 
\begin{equation}\label{eq:intequality}
    \partial_{t}(\rhoa \ea )+\nabla_{\bx}\cdot(\rhoa \ea \bua +\bq_h )+\textnormal{trace}(\Pressurea\nabla\bua) = - 2  \int  \limits_{{\mathcal S}(t)}\phiab (\bx,\by)\ea(t,\bx)\rho\rhobp \dy,
\end{equation}
with the so-called heat flux $\ds \bq_h:=\lim_{\Na  \rightarrow\infty}\frac{1}{2}\int |\bv-\bu|^2(\bv-\bu)\fa(t,\bx,\bv)\dv$.
Formally, any kinetic-based pressure tensor is in particular an entropic pressure, in the sense of satisfying the  \emph{equality} \eqref{eq:intequality}.
 But here one encounters the familiar problem of lack  closure which arises  whenever one is  dealing with  the highest truncated $\bv$-moments of $\fa$. In classical particle dynamics, the closure problem is resolved by  compatibility with  a preferred state of thermal equilibrium, a ``Maxwellian'' induced by the thermal equilibrium of the system, \cite{Lev1996,Gol1998,Cer2003,Vil2003}.
 In the current setup the agent-based dynamics, however, \eqref{eq:CS} governs \emph{active matter} made of `social particles' which admit no universal Maxwellian closure.  Then, there are multiple reasons which  led us to postulate the corresponding entropy \emph{inequality} \eqref{eq:meso-pressure}.
 
 \smallskip\noindent
 {\bf Scalar pressure}. We discuss the case of scalar pressure law $\ds {\Pressure}={\pressure}{\mathbb I}$. A large part of the existing literature on swarming \emph{assumes} a mono-kinetic closure,
  \begin{equation}\label{eq:mono-closure}
  \fa(t,\bx,\bv) \stackrel{N\rightarrow \infty}{\longrightarrow} \rho(t,\bx)\delta(\bv-\bu(t,\bx)),
  \end{equation}
   which is realized in terms of zero pressure, $\pressure=0$, e.g.,  \cite{HT2008,CFTV2010,FK2019,NP2021,Shv2021} and the references therein.  
 We mention the derivation  from  first principles \cite{Bia2012}, the isentropic closure, $\pressure=\rho^\gamma$, of \cite{KMT2013,KMT2015,KV2015,Cho2019,TCGW2020,Shv2022}, or equations of state fitted by observation that can be found in \cite{Sin2021} as  examples for  detailed thermodynamic closures for scalar pressure laws in the form of \emph{equality} in \eqref{eq:scalar-p} below.\newline
The notion of entropic pressure covers all these entropic examples of scalar pressure laws, as it applies to a broad class of pressure laws satisfying the entropy inequality postulated in \eqref{eq:meso-pressure}, but otherwise require  no algebraic closure. Indeed, our  notion of entropic pressure becomes more transparent in scalar case $\Pressure=\pressure{\mathbb I}$ where  the inequality postulated in \eqref{eq:meso-pressure}  for $\ds \pressure:=\frac{2}{d}\rho\ea$ reads (assuming  no heat flux $\bq=0$),
\begin{equation}\label{eq:scalar-p}
\partial_t  \pressure + \nabla_\bx\cdot (\pressure\bu)+\frac{2}{d}\pressure\nabla_\bx\cdot \bu \leq -2\pressure\int  \limits_{{\mathcal S}(t)}\phiab (\bx,\by)\rho(t,\bxp)\dy.
\end{equation}
Formal manipulation,  \eqref{eq:scalar-p}$ \times \rho^{-\gamma} - $ \eqref{eq:hydro}$\ds {}_1 \times\gamma \rho^{-\gamma-1}\pressure$ with $\gamma=1+\frac{2}{d}$, leads to the  equivalent entropic statement for $S=\ln\big(\pressure\rho^{-\gamma}\big)$,
\begin{equation}\label{eq:S}
\partial_t (\rho S) + \nabla_\bx\cdot (\rho \bu S) \leq -2\int  \limits_{{\mathcal S}(t)}\phiab (\bx,\by)\rho(t,\bx)\rho(t,\bxp)\dy, \qquad S:=\ln\big(\pressure\rho^{-(1+\frac{2}{d})}\big).
\end{equation}
We point out that the inequality \eqref{eq:S} is the \emph{reversed} entropy inequality encountered for $-S$ in compressible Euler equations. The difference,  which was already noted in \cite[\S6]{HT2008}, is due to different states of thermodynamic equilibria.

\noindent
{\bf Entropic energy dissipation}. An entropy inequality is intimately connected with the \emph{irreversibility} of the underlying process, e.g.,  the enlightening discussion in \cite[\S2.4]{Vil2003}. In the present context of hydrodynamics alignment, the entropy inequality  
\eqref{eq:meso-pressure}, or in its equivalent form \eqref{eq:entropy-ineq}, yields
\begin{equation}\label{eq:energy-dissipate}
 \begin{split}
 \ddt  \int \limits_{{\mathcal S}(t)}\rhoa \Ea\dx +  
 \int \limits_{\partial{\mathcal S}(t)}&\Big(\rhoa \Ea\bua\cdot {\mathbf n} +(\Pressurea\bua)\cdot{\mathbf n} +\bqa\cdot{\mathbf n}\Big){\d}S   \\
  & \leq 
 -  \iint\limits_{{\mathcal S}(t)\times{\mathcal S}(t)}\phiab (\bx,\by)\big(|\bua|^2-\bua\cdot\bubp +2\ea\big)\rhoa\rhobp \dx\dy\\
  & =- \hf\iint\limits_{{\mathcal S}(t)\times{\mathcal S}(t)}\phiab (\bx,\by)\Big(|\bubp-\bu|^2 +2\ea+2\ebp\Big)\rhoa\rhobp \dx\dy <0,
\end{split}
 \end{equation}
which reflects the dissipativity of the total energy $\ds \int \rhoa \Ea\dx$.
Thus, the entropy inequality \eqref{eq:meso-pressure} complements the 
balance laws in \eqref{eqs:hydro} to  govern the energy dissipation \eqref{eq:energy-dissipate}.
This is reminiscent of P.-L. Lions' notion of  dissipative solutions in the context of the Euler equations \cite[\S4.4]{Lio1996}. 
 
  One of the  main aspects of this work is dealing with arbitrary pressure, without  any specifics about the second-order closure for $\Pressurea$. 
The definition of entropic pressure in \eqref{eq:meso-pressure} is not concerned  with the detailed balance  of internal energy. Instead, its main purpose is to secure the dissipative nature of the total energy, $\displaystyle \rhoa \Ea$.  This partially echoes  Vicsek \& Zaferis who argued that in the context of collective motion ``\emph{The source of energy making the motion possible ... are not relevant}'' \cite[\S1.1]{VZ2012}. Here,  we abandon a closure in the form of thermal equality \eqref{eq:intequality} and instead, retain the inequality  postulated in  \eqref{eq:meso-pressure}, compatible with the dissipativity of internal fluctuations which we argued for in \cite[p. 501]{Tad2021}.   In particular, our definition of a  pressure in \eqref{eq:meso-pressure} can be realized in any intermediate scale between the microscopic agent-based description, \eqref{eq:CS}, and the macroscopic hydrodynamics \eqref{eqs:hydro}, and hence can be viewed as ``mesoscopic''. These considerations become even more pronounced when we extend our discussion to a larger class of so-called $p$-alignment hydrodynamics.
\section{$p$-alignment}\label{sec:p-alignment}
We  begin with the agent-based description,
\begin{equation}\label{eq:palignment}
    \left\{\begin{array}{c}
    \begin{split}
        \ddt\bx_{i}(t) &= \bv_{i}(t),   \\
        \ddt \bv_{i}(t) &= \frac{1}{\Nb}\sum_{j = 1}^{\Nb}\phiab_{ij}(t)|\bv_{j}(t)-\bv_{i}(t)|^{2p-2}(\bv_{j}(t)-\bv_{i}(t)),
    \end{split}      
    \end{array}\right. \quad i=1,2,\ldots \Na.
\end{equation}
 The case $p=1$ coincides with the Cucker-Smale model \eqref{eq:CS}, while for $p>1$, the alignment term  on the right of \eqref{eq:palignment} corresponds to weighted graph $2p$-Laplacian\footnote{To simplify computations we proceed  with $2p$-Laplacians rather than $p$-Laplacians.} which is found in recent applications of neural networks \cite{FZB2021}, spectral clustering \cite{BH2009}, semi-supervised learning \cite{ST2019},\cite{Fu2021}. In the context of alignment dynamics it was introduced  in \cite{HHK2010,CCH2014}.  We were motivated by the example of Elo rating system, \cite{JJ2015,DTW2019}, in which the alignment of scalar ratings $\{q_i\}$ is governed by odd function of local gradients $(q_j-q_i)$, e.g.,  $|q_j-q_i|^{2p-2}(q_j-q_i)$.\newline
  The long time behavior of the  $p$-alignment model  with $p>1$ is distinctly different from the `pure' alignment model when $p=1$ (and there is yet  a different behavior for $0\leq p <1$ which  we comment in remark \ref{rem:zalignment} below).
The large-crowd dynamics associated with \eqref{eq:palignment} is captured by the corresponding hydrodynamic description \myr{\myout{(again, its detailed derivation outlined in appendix \ref{sec:hydro-des} below)}}
 \begin{subequations}\label{eqs:phydro}
\begin{equation}\label{eq:phydro}
    \left\{\begin{array}{c}
    \begin{split}
    & \partial_{t}\rhoa+\nabla_\bx\cdot(\rhoa\bua) = 0,\\
    & \partial_{t}(\rhoa\bua)+\nabla_\bx\cdot(\rhoa\bua\otimes\bua+\Pressurea) = \aligna_p(\rho,\bub), 
    \end{split}      
    \end{array}\right.
\end{equation}
 with $p$-alignment term
 \begin{equation}\label{eq:palign}
  \aligna_p(\rho,\bub):=\int \limits_{{\mathcal S}(t)}\phi(\bx,\bxp)|\bubp-\bub|^{2p-2}(\bubp-\bu)\rho\rhobp\dxp, \qquad p\geq 1.
 \end{equation}
 \end{subequations}
 \myr{\begin{remark}[{\bf General $p$-alignment terms}]\label{rem:zero-average-Ap}
A detailed derivation of the $p$-alignment term $\aligna_p(\rho,\bu)$ in \eqref{eq:palign} is outlined in appendix \ref{sec:hydro-des}. This kinetic-based derivation is compatible with the mono-kinetic closure \eqref{eq:mono-closure}. In fact, our line of arguments below  does not require the detailed form of $\aligna_p(\rho,\bu)$, except for satisfying two `structural' conditions. The first condition requires that it has a zero average
 $ \ds \int \limits_{{\mathcal S}(t)} \aligna_p(\rho,\bu)(t,\bx)\dx=0$.
 This clearly holds for the $p$-alignment \eqref{eq:palign}, and in fact it holds for \emph{any} kinetic closure; see \eqref{eq:zero-kaverage-Ap} below. The second and essential condition  requires a $p$-alignment term which induces an entropic pressure. We discuss this notion of entropic pressure in context of $p$-alignment next.
 \end{remark}
 }
 
 \smallskip\noindent
 We assume that $\Pressure$ belongs to a class of entropic pressures, whose definition is adapted to the case of $p$-alignment.
 \begin{definition}[{\bf Entropic pressure for $p$-alignment}]\label{def:meso-ppressure}
We say that $\Pressurea$ is a entropic pressure  associated with \eqref{eqs:phydro} if it has a non-negative  trace, $\rhoa \ea:=\hf \text{trace}(\Pressurea)\geq0$, satisfying
\begin{equation}\label{eq:meso-ppressure}
    \partial_{t}(\rhoa \ea )+\nabla_{\bx}\cdot(\rhoa \ea \bua +\bqa )+\textnormal{trace}(\Pressurea\nabla\bua) \leq - \hf  \int  \limits_{{\mathcal S}(t)} \phiab(\bx,\by)\big((2\ea)^p+(2\ebp)^p\big)\rho\rhobp \dy.
\end{equation}
Here $\bqa$ is an arbitrary $C^1$-flux.
\end{definition}

Definition \ref{def:meso-ppressure} is motivated by the underlying kinetic formulation, where one encounters  the $p$-alignment quantity, see appendix \ref{sec:meso-ppressure} below\footnote{Here and below we abbreviate
$\square':=\square(t,\by,\bvp)$},
\[
-\frac{1}{2}\int  \limits_{{\mathcal S}(t)} \phi(\bx,\by)\iint \limits_{\R^d\times \R^d}|\bv-\bvp|^{2p}f_Nf'_N\dv\dvp\dxp.
\]
One cannot close the  kinetic expression $\ds \iint |\bv-\bvp|^{2p}f_Nf'_N\dv\dvp$ in terms of the thermodynamic quantity $\ea$, without taking into account a more detailed thermodynamic information i.e.,  higher moments of the empirical distribution $f_N$. It is here that we abandon the detailed thermal equality in favor of the inequality which follows from polarization, 
$\bv-\bvp\equiv (\bv-\bu)+(\bu-\bu')+(\bu'-\bvp)$,
\[
\begin{split}
-\frac{1}{2}&\iint |\bv-\bvp|^{2p}f_Nf'_N\dv\dvp \\
 & \leq -\frac{1}{2}\Big(\iint \Big(|\bv-\bua|^2+ |\bvp-\bu'|^2\Big)f_Nf'_N\dv\dvp\Big)^p \big(\rho\rho'\big)^{-\tfrac{p}{p'}}
 \stackrel{N\rightarrow \infty}{\longrightarrow} -\hf \big((2\ea)^p+(2\ebp)^p\big)\rho\rho'.
\end{split}
\] 
This leads to the corresponding  term  of $p$-entropic  pressure postulated on the right of \eqref{eq:meso-ppressure}.\newline
The special case of pure alignment, $p=1$, offers an alternative derivation where polarization implies the equality, consult \eqref{eq:casepeq1} below
\[
\begin{split}
\iint &(\bv-\bu)\cdot(\bv-\bvp)\fa \fbp \dvp\dv \\
   & =   -\iint |\bv-\bu|^2\fa \fbp\dvp\dv  - \int (\bv-\bu)\fa \dv \cdot \int (\bu-\bvp)\fbp\dvp  \stackrel{N\rightarrow \infty}{\longrightarrow} -2\ea\rho\rho',
   \end{split}
\]
which in turn formally yields  the entropy  \emph{equality} \eqref{eq:intequality},
 \begin{equation}\label{eq:meso-pressureb}
    \partial_{t}(\rhoa \ea )+\nabla_{\bx}\cdot(\rhoa \ea \bua +\bqa )+\textnormal{trace}(\Pressurea\nabla\bua) = - 2\int  \limits_{{\mathcal S}(t)} \phiab (\bx,\bxp)\ea\rho\rhobp \dxp.
\end{equation}
Thus, while  for $p=1$ the inequality of entropic pressure \eqref{eq:meso-pressure} could be viewed as a matter of choice made in the equalities \eqref{eq:intequality} or \eqref{eq:meso-pressureb},  for  $p>1$ the entropic inequality \eqref{eq:meso-ppressure} is a necessity in order to have a macroscopic interpretation of an entropic pressure.

\begin{remark}[{\bf Local vs. global flux}]\label{rem:localvsglobal}
We observe that the entropic statement for $p$-alignment  \eqref{eq:meso-ppressure} with $p=1$  is a  symmetric version of the entrropic inequality of `pure' alignment,  \eqref{eq:meso-pressure}. Apparently,  the two definitions do not agree when $p=1$, but in fact,  their difference is encoded in different fluxes $\bq$. In particular, while the entropic pressure  in pure alignment \eqref{eq:meso-pressure} is encoded in terms of a \emph{local} heat flux, $\bq_h$ in \eqref{eq:heat-flux} below, the case of $p$-alignment \eqref{eq:meso-ppressure} requires a \emph{global} flux, $\bq_h+\bq_\phi$ in \eqref{eq:penergy-balance} below. Alternatively, we  could be less `pedantic'  and combine both cases of alignment and of $p$-alignment under the same  notion of  entropic pressure inequality
\[
    \partial_{t}(\rhoa \ea )+\nabla_{\bx}\cdot(\rhoa \ea \bua +\bqa )+\textnormal{trace}(\Pressurea\nabla\bua) \leq - 2^{p-1} \int  \limits_{{\mathcal S}(t)} \phiab(\bx,\by)\ea^p\rho\rhobp \dy, \qquad p\geq 1.
\]
This will not affect any of the follow-up results. 
\end{remark}

\noindent
Of course,  a general $C^1$ flux, $\bq$, can also `absorb' the convective term $\rho\ea\bu$; our main  focus is in the global dissipative structure entailed by \eqref{eq:meso-ppressure}.

\medskip\noindent
{\bf Entropic energy dissipation in $p$-alignment}.  Following the same formal manipulations  as before for $p=1$, see \eqref{eq:energy-kin}, yield
 \[
 \partial_t \Big(\frac{\rhoa}{2}|\bua|^2\Big) + \nabla_\bx\cdot \Big(\frac{ \rhoa}{2}|\bua|^2\bua +\Pressurea\bua\Big) - \textnormal{trace} \big(\Pressurea \nabla\bua\big)  \leq 
 -  \int \limits_{{\mathcal S}(t)}\phiab (\bx,\by)|\bu-\bubp|^{2p-2}\bua\cdot(\bua-\bubp)\rho\rhobp \dy.
\]
 Adding \eqref{eq:meso-ppressure} and integrating we find
\begin{equation}\label{eq:penergyinq}
 \begin{split}
 \ddt \int \limits_{{\mathcal S}(t)}&\rho E(t,\bx)\dx  + \int \limits_{\partial{\mathcal S}(t)}\Big(\rho E\bu\cdot \na + (\Pressure\bu)\cdot\na+\bq\cdot\na\Big){\d}S \\
  & \leq - \iint\limits_{{\mathcal S}(t)\times{\mathcal S}(t)}\phiab (\bx,\by)\Big(|\bu-\bubp|^{2p-2}\big(|\bua|^2-\bua\cdot\bubp) +\hf\big((2\ea)^p+(2\ebp)^p\big)\Big)\rho\rhobp \dy\\
  & = -\hf\iint \limits_{{\mathcal S}(t)\times{\mathcal S}(t)}\phi(\bx,\bxp)\Big(|\bu'-\bu|^{2p} +(2\ea)^p+(2\ebp)^p\Big)\rho\rho'\dx\dxp <0.
  \end{split}
 \end{equation}
which extends  the dissipativity statement of `pure' alignment  in the case $p=1$ in \eqref{eq:energy-dissipate}.

\ifx
\begin{remark}[{\bf Thermodynamic equilibrium}]
  Observe that the notion of entropic pressure is invariant under addition of rank-one matrices. Indeed, fix a unit vector $\bw\in \R^d$ and set $\ds \Pressure_\eta = \Pressure +\frac{2\eta}{d}\bw\bw^\top, \ \eta>0$ then the mass equation \eqref{eq:hydro}${}_1$ implies that \eqref{eq:meso-pressure} holds with $\rho e_\eta :=\rho\ea+\eta$ and $\bq_\eta:=\bq-\eta\bu$,
  \[
  \begin{split}
     \partial_{t}(\rhoa e_\eta ) +&\nabla_{\bx}\cdot(\rhoa e_\eta \bua +\bqa_\eta)  +\textnormal{trace}(\Pressure_\eta\nabla\bua) \\
      & =  \partial_{t}(\rhoa \ea ) +\nabla_{\bx}\cdot(\rhoa \ea \bua +\bqa )+\textnormal{trace}(\Pressure\nabla\bua) \\
      & \leq - 2^{2p-1}\rhoa \ea^p  \int \limits_{\Omegae} \phiab (\bx,\by)\rhobp \dy < - 2^{2p-1}\rhoa e^p_\eta  \int \limits_{\Omegae} \phiab (\bx,\by)\rhobp \dy, \quad \eta>0. 
     \end{split}
  \]
    \end{remark}
\fi
 \section{Swarming}\label{sec:swarming} 
 The hydrodynamic alignment \eqref{eqs:hydro}  occupies a distinct `patch' of mass, 
\[
\Omegaa(t)=\textnormal{supp}\, \rhoa(t,\cdot).
\] 
We shall refer to this patch of mass  simply as a `crowd'  --- a continuum of agents which encodes the large-crowd dynamics associated with \eqref{eq:CS}.
In most of the existing literature on collective dynamics, the edge of such crowd is assumed to be `tailored'  to the surrounding vacuum so that ${\rhoa}(t,\cdot)_{{}|_{\partial\Omegaa}}=0$. Instead, we argue here  for a more realistic scenario in which the density inside the crowd remains strictly bounded away from vacuum, 
\begin{equation}\label{eq:vaccum}
\min_{\bx\in\Omegaa(t)}\rhoa(t,\bx)\geq \rho_->0,
\end{equation}
while  its  boundary 
forms a shock discontinuity, moving with velocity ${\bua}_{|\partial \Omegaa}$.  A detailed discussion on the nature of boundary conditions (BCs) for such crowds is beyond the scope of this work (see \cite{AC2021a} for the special  one-dimensional case with $p=\rho$). Instead, we argue \eqref{eqs:hydro} augmented  with Neumann BCs 
\begin{equation}\label{eq:BCs}
\bua\cdot  \na_{|\partial\Omegaa}=0, \quad \Pressurea\na_{|\partial\Omegaa}=0, \quad \textnormal{and} \quad \bqa\cdot\na_{|\partial\Omegaa}=0.
\end{equation}
In particular, it follows  that the total mass of the crowd, $M=M(t)$, is conserved in time,
\begin{equation}\label{eq:conserve-mass}
\Ma(t) :=\int \limits_{\Omegaa(t)}\rhoa(t,\bx)\dx \equiv \Ma_0, 
\end{equation}
and by  the symmetry of $\phiab (\cdot,\cdot)$
\[
\ddt  \int \limits_{\Omegaa(t)}\rhoa\bua\dx=-\int \limits_{\partial\Omegaa(t)} \left(\bua\otimes \bua\cdot \na +\Pressurea\na\right)\rho {\d}S \ \ - \!\!  \iint \limits_{\Omegaa(t)\times\Omegab(t)} \phiab (\bx,\by)|\bubp-\bu|^{2p-2}(\bubp -\bua)\rhoa\rhob'\dx\dxp =0,
\] 
and hence the total momentum of the crowd, $\bma=\bma(t)$, is also conserved,\footnote{\myr{This is the only stage which requires the zero-average  $p$-alignment term argued in remark \ref{rem:zero-average-Ap}, 
\[
 \int \limits_{{\mathcal S}(t)} \aligna_p(\rho,\bu)(t,\bx)\dx =
 \iint \limits_{\Omegaa(t)\times\Omegab(t)} \phiab (\bx,\by)|\bubp-\bu|^{2p-2}(\bubp -\bua)\rhoa\rhob'\dx\dxp =0,
 \]
 which in turn implies conservation of total momentum $\bma(t)=\bma_0$.
}}    
\begin{equation}\label{eq:conserve-momentum}
 \bma(t):=\int \limits_{\Omegaa(t)} \rhoa(t,\bx)\bua(t,\bx)\dx \equiv \bma_0.
\end{equation}
 Finally,  \eqref{eq:penergyinq} yields that the  total energy is non-increasing
 \begin{equation}\label{eq:total-energy}
 \begin{split}
 \ddt \int \limits_{{\mathcal S}(t)}\rho E(t,\bx)\dx  \leq -\hf\iint \limits_{{\mathcal S}(t)\times{\mathcal S}(t)}\phi(\bx,\bxp)\Big(|\bu'-\bu|^{2p} +(2\ea)^p+(2\ebp)^p\Big)\rho\rho'\dx\dxp.
  \end{split}
 \end{equation}
 In particular, we have  the space-time  \emph{enstrophy} bound 
 \begin{equation}\label{eq:ens}
 \int \limits_0^t \iint \limits_{{\mathcal S}(t)\times {\mathcal S}(t)} \phi(\bx,\bxp)\Big(|\bubp-\bu|^{2p} +(2\ea)^p+(2\ebp)^p\Big)\rho\rho' \dx\dxp{\d}t \leq C^2_0:=2\int  \limits_{{\mathcal S}(0)}\rho_0 E_0\dx.
 \end{equation}
 
 \medskip\noindent
 {\bf Flocking}. A characteristic feature of alignment dynamics is the emergence of coherent structure with limiting velocity  $\bu_\infty$ such that
  \begin{equation}\label{eq:u-convergence}
 \bu(t,\bx)-\bu_\infty(t,\bx) \stackrel{t\rightarrow \infty}{\longrightarrow}0,
 \end{equation}
 and the corresponding limiting density, $\rho_\infty$. This is typical in \emph{flocking} phenomena.
  In the  present context of the hydrodynamic alignment \eqref{eqs:hydro},  the limiting behavior of the dynamics  \eqref{eqs:hydro} can only approach the time-invariant mean velocity $\ds \bu_\infty=\bubar:=\frac{\bm_0}{M_0}$ with  a limiting density  carried out as a traveling wave $\rho_\infty(\bx-\bubar t)$, \cite[\S2]{ST2017b}. The presence of additional repulsion, attraction and external forces introduce a `richer' set  of  possible emerging limiting  configurations, e.g., \cite{CDMBC2007}; for example, alignment with quadratic forcing approaching an harmonic oscillator $\ddot{\bu}_\infty(t)+a^2\bu_\infty(t)=0$, \cite[\S2.4]{ST2020a}.
 The precise notion of flocking convergence in \eqref{eq:u-convergence} may vary. Ideally, we seek uniform convergence. A more relaxed notion  of $L^2_\rho$-convergence becomes accessible by studying \emph{energy fluctuations}, see the next section \ref{sec:fluctuations},
\[
  \int \limits_{{\mathcal S}(t)} |\bu-\bu_\infty|^2\rho\dx \stackrel{t \rightarrow \infty}{\longrightarrow}0.
\] 
 In practice, as we shall see below, the analysis may gain by a combination of the two.\newline  
 The limiting configuration is supported on 
$ {\mathcal S}_\infty(t):=\textnormal{supp}\, \rho_\infty(t,\cdot)$.
 For example, ${\mathcal S}_\infty(t)$ is a Dirac mass in presence additional attractive forces \cite[Theorem 1]{ST2021}.
Ideally, we are interested to trace the shape of the boundary $\partial{\mathcal S}_\infty(t)$, but this seems to be out of reach in the current literature (but see \cite{LLST2022}). In general, one expects that  alignment is at least strong enough to keep the dynamics  contained in a finite ball, 
 \[
  D(t)\leq D_+<\infty, \qquad D(t):=\max_{\bx,\bxp \in {\mathcal S}(t)}|\bx-\bxp|.
  \]
In practice we may need to  address to a more accessible notion of diameter which allows a slow time growth, $D(t) \leq C_D\myangle{t}^{\gamma}$ with some fixed $\gamma>0$. 

\section{Decay of energy fluctuations}\label{sec:fluctuations}
We study the hydrodynamics of  the $p$-alignment,  \eqref{eqs:phydro}, assuming it admits a strong entropic solution, \eqref{eq:meso-ppressure}; see further comments on \ref{H1} in section \ref{sec:ens-and-dispersion} below.\newline
Consider  the \emph{energy fluctuations}, \cite[\S5]{HT2008}, \cite{Tad2021},
\begin{equation}\label{eq:eflock}
\delE(t):= \frac{1}{2M} \iint \limits_{\Omegaa(t)\times \Omegaa(t)}\Big(\hf |\bua (t,\bx)-\bu(t,\bxp)|^2+\ea(t,\bx)+\eb(t,\bxp)\Big)\rhoa(t,\bx)\rhob(t,\bxp)\dx\dxp.
\end{equation}
It can be expressed in the equivalent form,\footnote{Specifically
\[ 
\begin{split}
\frac{1}{M}\iint  \limits_{\Omegaa(t) \times \Omegab(t)} &\hf|\bua(t,\bx)-\bub(t,\bxp)|^2\rhoa\rhobp \dx\dy \\
& =\frac{1}{M}\iint \limits_{\Omegaa(t) \times \Omegab(t)} \left(\hf|\bua(t,\bx)-\bubar|^2 + (\bua-\bubar)\cdot(\bubar-\bubp )
+\hf|\bu (t,\bxp)-\bubar|^2\right)\rhoa\rhobp \dx\dy \\
 &= \int \limits_{{\mathcal S}(t)} |\bua(t,\bx)-\bubar|^2\rhoa(t,\bx)\dx. 
\end{split}
\]
}
$\ds \delE(t)=  \int \limits_{\Omegaa(t)}\Big(\hf|\bua (t,\bx)-\bubar(t)|^2+\ea(t,\bx)\Big)\rhoa(t,\bx)\dx$.
Thus, $\delE(t)$ reflects macroscopic velocity fluctuations  $\ds \int\limits_{\Omegaa(t)} \hf|\bu-\bubar(t)|^2\rho(t,\bx)\dx$ around the mean velocity, $\ds \bubar(t):=\frac{1}{M}\int\limits_{\Omegaa(t)} \rho\bu(t,\cdot)\dx = \frac{\bm}{M}$, and\myr{\mycancel{,}} in the context of kinetic formulation \eqref{eq:pressure}--\eqref{eq:internale}, it also reflects the microscopic  velocity fluctuations, 
$\ds \rhoa\ea  = \lim_{N\rightarrow \infty} \int \hf|\bv-\bu|^2\fa(t,\bx,\bv)\dv$.
We have the following decay bound on energy fluctuations 
\begin{equation}\label{eq:phydroRic}
\ddt \delE(t)  \leq  -2^pM^{2-p}k(D(t))\big(\delE(t)\big)^p. 
\end{equation}
The derivation follows  the energy inequality \eqref{eq:total-energy}. Noting that $\ds  \delE(t)\equiv \int \limits_{{\mathcal S}(t)} \rhoa\Ea \dx-\frac{1}{2M}\big| \bma\big|^2$ with total mass and total momentum which are  conserved in time, $\Ma(t)= \Ma_0, \  \bma(t)=\bma_0$,  we end up with,
\begin{equation}\label{eq:delEpdecay}
\begin{split}
\ddt &\delE(t)   = \ddt \int \limits_{{\mathcal S}(t)}\rhoa \Ea(t,\bx)\dx\\
 & \leq  -\hf\iint \limits_{{\mathcal S}(t)\times {\mathcal S}(t)} \phi(\bx,\bxp)\Big(|\bu-\bubp|^{2p}+(2\ea)^p+(2\ebp)^p\Big)\rho \rhobp\dx\dxp\\
& \leq -\iint \limits_{{\mathcal S}(t)\times {\mathcal S}(t)} \phi(\bx,\bxp)\Big(\frac{1}{2}|\bu-\bubp|^2+\ea(t,\bx)+\eb(t,\bxp)\Big)^p\rho \rhobp\dx\dxp \\
 & \leq -k(D(t))\Big(\iint \limits_{{\mathcal S}(t)\times {\mathcal S}(t)} \hspace*{-0.2cm}\Big(\frac{1}{2}|\bu-\bubp|^2+\ea+\ebp\Big)\rho\rho' \dx\dxp\Big)^p \Big(\iint \limits_{{\mathcal S}(t)\times {\mathcal S}(t)} \hspace*{-0.2cm}\rho\rho' \dx\dxp\Big)^{-\tfrac{p}{p'}} \\
 &=  - 2^pM^{2-p}k(D(t))\big(\delE(t)\big)^p.
\end{split}
\end{equation}
The first inequality on the right quotes \eqref{eq:total-energy}; the second follows from Jensen inequality and the third from  H\"{o}lder inequality and the obvious radial bound \eqref{eq:radial}, $\phi(\bx,\by)\geq k(D(t))$. Integration of \eqref{eq:delEpdecay} yields the following.

\begin{theorem}\label{thm:main1} Let $(\rho,\bu,\Pressure)$   be a strong  solution\footnote{That is, $\big(\rhoa(t,\cdot),\bua (t,\cdot),\Pressurea (t,\cdot)\big)$ has sufficient smoothness --- say  $\in L_+^{\infty}\cap L^{1})\big({\mathcal S}(t)\big)\times W^{1,\infty}\big({\mathcal S}(t)\big)\times W^{1,\infty}\big({\mathcal S}(t)\big)$, so that \eqref{eqs:hydro} can be interpreted in a  pointwise sense.} of the hydrodynamic $p$-alignment  \eqref{eqs:phydro}, satisfying the entropy condition \eqref{eq:meso-ppressure},  and subject to compactly supported initial data, $(\rhoa_0,\bua_0,\Pressurea_0)$ with $D_0<\infty$, and boundary conditions \eqref{eq:BCs}.
Then the  \emph{energy fluctuations} $\delE(t)$
admits the bound
\begin{equation}\label{eq:decayfluc}
\begin{split}
     \delE(t) \leq \left\{\begin{array}{ll}
     \ds exp\left\{-2M\int \limits_{0}^{t}\kab\big(D(s)\big)\d{s}\right\}\delE(0), & p=1\\ \\
    \ds  \frac{1}{\ds \left\{(p-1)2^pM^{2-p}\int \limits_0^t k(D(s))\d{s}\right\}^{\tfrac{1}{p-1}}}, & p>1.
     \end{array}\right.
\end{split}
\end{equation}
\end{theorem}
\noindent
The result applies to $p$-alignment dynamics with general class of entropic pressure tensors satisfying  \eqref{eq:meso-ppressure} (noting that \eqref{eq:meso-pressure} for $p=1$ yields the same energy decay  \eqref{eq:penergyinq}). We refer to such solutions as 
`entropic solutions'. The  symmetric communication protocol $\phiab$ in  \eqref{eq:radial}  need not be metric nor bounded and no assumption of a uniform velocity bound is made.

We close by noting that the bound \eqref{eq:decayfluc}${}_2$ depends on the initial mass $M$ but otherwise it is independent of the initial fluctuations $\delE(0)$ --- a typical scenario for the Ricatti's type inequality \eqref{eq:phydroRic} with $p>1$.
\subsection{Heavy-tailed kernels}\label{sec:heavy-tailed}
The bound \eqref {eq:decayfluc} reflects  a competition between the expansion rate of the diameter of the crowd,
$D(t)$, and the decay rate in   its communication strength, $\kab(r)$:  their composition is required to have a  non-integrable ``heavy-tail'' in order to enforce   $L^2_\rho$-flocking decay. 
We make these considerations precise in our next statement.\newline
{\bf Communication kernels of order $\beta\geq 0$}. There exist  constants $C_k>0, R>0$ such that
\begin{equation}\label{eq:pconnect}
  \phiab (\bx,\by) \geq  \kab (|\bx-\by|)\ \ \textnormal{with} \ \ \left\{\begin{array}{l}
   \ds \int \limits_{|\bx|\leq R} k(|\bx|)\dx  < \infty, \\ \\
  \kab(r) = C_k \myangle{r}^{-\beta},   \quad r\geq R.
   \end{array}\right.
\end{equation}
This emphasizes the fact that besides the mere requirement for intractability of  $\phi$  near the origin --- only its tail behavior  matters.\newline
{\bf Notations}. We use the following two constants. 
We let $C_R$  
denote a constant, with different values in different contexts, depending of $R$ as well as on the other fixed parameters $\beta,\gamma,...$ and possibly $p>1$. Also, we denote the `scaled mass'
\[
M_p:=\left\{\begin{array}{ll}
2MC_kC_D^{-\beta}, & p=1,\\ \\
\big(2^pM^{2-p}C_kC_D^{-\beta}\big)^{-\tfrac{1}{p-1}}, & p> 1.
\end{array}\right.
\]
\begin{corollary}[{\bf Decay of $L^2_\rho$-energy fluctuations}]\label{cor:1}
Let $(\rhoa,\bua,\Pressure)$   be a strong entropy solution of the hydrodynamic $p$-alignment system \eqref{eqs:phydro},\eqref{eq:meso-ppressure}, $p\geq 1$,  with communication kernel  $\phiab$ of order $\beta\geq 0$, 
\eqref{eq:pconnect}.
Assume that the ``crowd'' disperses at a rate of order $\gamma\geq 0$,
\begin{equation}\label{eq:sizeofDt}
D(t) \leq C_D(1+t)^\gamma, \qquad \gamma\geq 0, \qquad D(t)=\max\{|\bx-\bxp|, \ \bx,\bxp\in \textnormal{supp}\, \rho(t,\cdot)\}.
\end{equation}
 If the heavy-tail condition holds in the sense that  $\beta\gamma < 1$,
  then there is long time flocking behavior 
  such that the following  decay  bound holds
\begin{equation}\label{eq:flocking_decay}
\begin{split}
\delE(t) 
  \leq \left\{\begin{array}{ll}
  \ds C_R\, exp\big\{- M_1 t^{(1-\beta\gamma)}\big\}\delE(0), & p=1\\ \\
  \ds  C_R M_p\, t^{-\tfrac{1-\beta\gamma}{p-1}}, & p>1.
  \end{array}\right.
  \end{split}
\end{equation}
\end{corollary}
\noindent
In case of pure alignment, $p=1$, \eqref{eq:flocking_decay}${}_1$ recovers an exponential  decay of fractional  order $1-\beta\gamma$, \cite[Corollary 1]{Tad2021}, while for $p>1$, \eqref{eq:flocking_decay}${}_2$ implies a Pareto-type decay of fractional order $\ds \frac{1-\beta\gamma}{p-1}$.
Thus, corollary \ref{cor:1} implies that for heavy-tailed kernels such that $\beta\gamma<1$, both  the macroscopic and microscopic fluctuations around the mean $\bubar(t)=\bubar_0$ decay to zero. In particular, this shows the \emph{trend towards equilibrium} of a kinetic-based hydrodynamics, as it decays towards mono-kinetic closure \eqref{eq:mono-closure}
\[
 \hf\int \limits_{{\mathcal S}(t)} \|\Pressure(t,\bx \|_2\dx = \int \limits_{{\mathcal S}(t)} \ea(t,\bx)\rho(t,\bx)\dx  \stackrel{t\rightarrow \infty}{\longrightarrow}0.
\]
A key aspect, therefore,  is to study the  possible expansion  of the  spatial diameter with time growth of order $\gamma$ (possibly depending on $\beta$), so that $\beta\gamma<1$.
This will occupy us in the rest of the work.  

\begin{remark} One can refine the statement of corollary \ref{cor:1} to include the borderline case, $\beta\gamma=1$ 
\end{remark}

\section{Flocking with mono-kinetic (``pressure-less'') closure}\label{sec:mono-kinetic}
One strategy for verifying flocking is to seek a  uniform bound on velocity, $u_+:= \max|\bu(t,\cdot)| <\infty$, which in turn implies a dispersion bound on the diameter of order $\lesssim (1+t)$,
\begin{equation}\label{eq:mono-dispersion}
\ddt D(t) \leq \delbu(t), \quad \delbu(t):=\max_{\bx,\bxp\in {\mathcal S}(t)}|\bu(t,\bxp)-\bu(t,\bx)| \ \  \leadsto \ \ D(t)\leq D_0+2u_+ t, 
\end{equation}
and then appeal to corollary \ref{cor:1} with $\gamma=1$. An instructive example for this line of argument is found in the prototype  case of  \emph{mono-kinetic closure},    $\Pressure =0$, 
 \begin{equation}\label{eq:hydro-mono-kinetic}
  \partial_{t}(\rhoa\bua)+\nabla_\bx\cdot(\rhoa\bua\otimes\bua) = \align_p(\rho,\bu).
 \end{equation}
   A main feature of the mono-kinetic closure  is that the resulting system    \eqref{eq:hydro-mono-kinetic} decouples into scalar transport equations,  
   \[
   u_t +\bu\cdot\nabla_\bx u = \int \limits_{{\mathcal S}(t)}\phi(\bx,\bxp)|\bubp-\bu|^{2p-2}(u'-u)\rhobp\dxp, 
   \]
in which case, the coercivity of the (scalar) $p$-alignment term on the right implies a maximum principle, $\max|\bu(t,\cdot)|\leq \max|\bu_0|$, hence
\[
D(t)\leq D_0+2u_+\!\cdot t, \qquad  u_+:= \max|\bu_0|.
\]
Appealing to corollary \ref{cor:1} with $\gamma=1$  implies that for heavy-tailed $\phi$'s of order $\beta<1$, there exists $C_R=C(R,D_0,u_+,\beta,p)$ such that
\[
\begin{split}
\delE(t) 
  \leq \left\{\begin{array}{ll}
  \ds C_R\, exp\big\{-2M(D_0+2u_+\!\cdot t)^{(1-\beta)}\big\}\delE(0), & p=1\\ \\
  \ds  C_RM_p\big(D_0+2u_+\!\cdot t\big)^{-\tfrac{1-\beta}{p-1}}, & p>1.
  \end{array}\right.
  \end{split}
\]
In fact, more is true --- a refined argument shows that for such heavy-tailed $\phi$'s of order $\beta<1$, the pressureless diameter remains uniformly bounded, $D(t) \leq D_+$, and hence corollary \ref{cor:1} applies with $\gamma=0$. To this end, we split out discussion, distinguishing between the  case of `pure' alignment, $p=1$, and the case of $p$-alignment $p>1$.
 \subsection{Flocking with pure alignment ($p=1$)}
We begin with the following pointwise bound of velocity fluctuations  which is reproduced in section \ref{sec:pointwise_mono-kinetic} below,
\begin{equation}\label{eq:delV-contarct}
\ddt \delbu(t) \leq -k(D(t))M\delbu(t), \qquad  \delbu(t)=\sup_{\bx,\by\in{\mathcal S}(t)}|\bua(\bx,t)-\bub(\by,t)|.
\end{equation}
In particular, $\delbu(t)\leq \delbu_0$ and hence  \eqref{eq:sizeofDt} holds with $\gamma=1$ in view of $D(t) \leq D_0 +\delbu_0\cdot t$.
Consequently,  for $\beta$-tailed kernels of order $\beta <  1$,   \eqref{eq:pconnect}, there exists a constant $c_R$ such that
\[
\displaystyle \int \limits_{\myr{0}}^{\myr{t}} k(D(s)){\d}s \geq \int \limits_{\myr{0}}^{\myr{R}} k(D(s)){\d}s + \int \limits_{\myr{R}}^{\myr{\max\{R,t\}}} k(D(s)){\d}s \geq \frac{c_R}{\myr{(1-\beta)\delbu_0}}(1+\delbu_0\cdot t)^{1-\beta}, \quad 0\leq \beta<1. 
\]
Revisiting \eqref{eq:delV-contarct} again  yields a \emph{decay}  of pointwise velocity fluctuations of fractional exponential order, $\ds \delbu(t)\leq \delbu_0 \,exp\{-c'_R (1+\delbu_0\cdot t)^{1-\beta}\}$ with $\ds c'_R=\frac{M}{(1-\beta)\delbu_0}c_R$, 
which in turn implies that the diameter remains uniformly bounded 
\[
\ddt D(t)\leq \delbu_0 \,e^{-c'_R (1+\delbu_0\cdot t)^{1-\beta}}\ \ \leadsto \ \ D(t) \leq D_+:= \delbu_0 \int_0^\infty e^{-c'_R (1+\delbu_0\cdot t)^{1-\beta}}{\d}t <\infty.
\] 
Alternatively, one can use the deceasing Liapunov functional of \cite{HL2009},
$\ds \delbu(t)+M\int_{D_0}^{D(t)}k(s){\d}s$ to conclude that any heavy-tailed kernel in the sense that $\ds \int k(s){\d}s =\infty$ implies $D(t)\leq D_+<\infty$.
 Thus, whenever $\beta <1$, then 
 corollary  \ref{cor:1} applies with $\gamma=0$  and $C_D=D_+$ and one   recovers the exponential decay of  mono-kinetic dynamics, \cite{CS2007a,HT2008,HL2009,CFTV2010,Shv2021}. 
 \begin{proposition}[{\bf Flocking for mono-kinetic alignment, $p=1$}]\label{prop:1}
Let $(\rhoa,\bua)$   be a strong solution of the  mono-kinetic alignment system \eqref{eqs:hydro} with
 ``heavy-tailed'' communication kernel  $\phiab$ of order $0\leq \beta< 1$, \eqref{eq:pconnect}.
There is long time flocking behavior with  decay rate 
\begin{equation}\label{eq:delE-revisit}
\int \limits_{\Omegaa(t)} |\bua(t,\bx)-\bubar|^2\rhoa(t,\bx)\d\bx \leq C_R \, e^{-M_1 t}\!\int \limits_{\Omegaa(t)} |\bua_0(\bx)-\bubar|^2\rhoa_0(\bx)\d\bx, \quad M_1=2MC_kD_+^{-\beta}.
\end{equation}
\end{proposition}
Moreover, integration of \eqref{eq:delV-contarct} implies pointwise bound on  the decay of velocity fluctuation
\begin{equation}\label{eq:vfluc-revisit}
 \max_\bx|\bua(t,\bx)-\bubar| \leq C_R \, e^{-M_1 t}\max_\bx  |\bua_0(\bx)-\bubar|.
\end{equation}
 \subsection{Flocking with $p$-alignment ($p>1$)}
Our starting point is  the pointwise bound of velocity fluctuations  corresponding to \eqref{eq:delV-contarct}, which is outlined in appendix \ref{sec:pointwise_pmono-kinetic}
\begin{equation}\label{eq:delV-pcontarct}
\ddt \delbu(t) \leq -\hf M k(D(t))(\delbu(t))^{2p-1}, \qquad  \delbu(t)=\sup_{\bx\in{\mathcal S}(t)}|\bua(\bx,t)-\bubar|.
\end{equation}
In particular, $\delbu(t)\leq \delbu_0$ implies $D(t)\leq D_0+2\delbu_0\cdot t$, that is, \eqref{eq:sizeofDt} holds with $\gamma=1$,
\[
D(t)\leq C_D(1+t), \qquad C_D=\max\{D_0,2\delbu_0\},
\]
\ifx
We now recall that $k(r)\geq C_k(1+r)^{-\beta}$ for $r>R$. Hence, using \eqref{eq:delV-pcontarct} for $t \gtrsim R$ where $k(D(t))\geq C_RAB^{-\beta}(1+t)^{-\beta}$  leads to 
\[
\ddt \big(\delbu(t)\big)^{2-2p} \geq C_RM_1(1+t)^{-\beta}, \quad  p>1,
\]
where, as before, $M_1=2MC_kC_D^{-\beta}$.
Integrating, we conclude with the flocking bound with decay rate of order $\frac{1-\beta}{2p-2}$.
\[
\delbu(t) \leq C_R \frac{1}{\big\{M_1(1+t)^{1-\beta}+(\delbu_0)^{2-2p}\big\}^{\frac{1}{2p-2}}} \leq C_R M_1^{-\frac{1}{2p-2}}(1+t)^{-\frac{1-\beta}{2p-2}}.
\]
\fi
and corollary \ref{cor:1} implies  $L^2_\rho$-decay rate of order $\ds \frac{1-\beta}{p-1}$. 
\begin{proposition}[{\bf Flocking for mono-kinetic alignment, $p>1$}]\label{prop:2}
Let $(\rhoa,\bua)$   be a strong solution of the mono-kinetic $p$-alignment system \eqref{eqs:phydro} with   ``heavy-tailed'' communication kernel  $\phiab$ of order $0\leq \beta< 1$, \eqref{eq:pconnect}.
Then there is long time flocking behavior with  decay rate 
\begin{equation}\label{eq:pflocking_decay}
\begin{split}
\delE(t) \leq 
   \ds  C_RM_p (1+t)^{-\frac{1-\beta}{p-1}}, \qquad  p>1, \quad 0\leq \beta<1.
    \end{split}
\end{equation}
\end{proposition}
\noindent
We can improve these bounds, at least in the restricted range $1<p<\nicefrac{3}{2}$.
To this end,  use an iterative argument starting with the $\gamma$-bound
\[
D(t) \leq C_D(1+t)^\gamma.
\]
Integrating \eqref{eq:delV-pcontarct} for $t \gtrsim R^{1/\gamma}$ where $k(D(t))\geq C_RC_kC_D^{-\beta}(1+t)^{-\beta\gamma}$  leads to 
\[
\ddt \big(\delbu(t)\big)^{2-2p} \geq C_R(p-1)M_1(1+t)^{-\beta\gamma},  \quad p>1, 
\]
where, as before, $ M_1=MC_kC_D^{-\beta}$.
We conclude with the flocking bound 
\[
\delbu(t) \leq C_R \frac{1}{\big\{M_1(1+t)^{1-\beta\gamma}+(\delbu_0)^{2-2p}\big\}^{\tfrac{1}{2p-2}}} \leq C_R M_1^{-\tfrac{1}{2p-2}}(1+t)^{-\tfrac{1-\beta\gamma}{2p-2}},
\]
and hence
\begin{equation}\label{eq:Dgamma}
\ddt D(t) \leq 2\delbu(t) \ \ \leadsto \ \ D(t)\leq D_0+ C_R2\frac{2M_1^{-\tfrac{1}{2p-2}}}{\gamma'}(1+t)^{\gamma'}, \quad \gamma':=\frac{2p-3}{2p-2}+\frac{\beta\gamma}{2p-2}.
\end{equation}
We distinguish between two cases. If $2p+\beta<3$ then after one iteration, starting with $\gamma=1$, we obtain
\[
\gamma'=\frac{2p-3+\beta}{2p-2} <0.
\] 
 If, however, $2p+\beta\geq3$ and $\beta<\nicefrac{1}{2}$ then $\frac{\beta}{2p-2}<1$ and, hence the  fixed point iterations $\gamma\mapsto \gamma'$ form a contraction, approaching  the negative value
 \[
 \gamma_\infty = \frac{2p-3}{2p-2-\beta}< 0, \qquad p<\nicefrac{3}{2}.
 \]
 In either case, the range $1 \myr{<} p < \nicefrac{3}{2}$ and $\beta<\nicefrac{1}{2}$   implies   that after finitely many iterations,  \eqref{eq:Dgamma} holds with $\gamma <0$ and we conclude that the diameter $D(t)$ remains uniformly bounded in time, $D(t)\leq D_+$, that is, \eqref{eq:sizeofDt} holds with $\gamma=0$ and $C_D=D_+$. Corollary \ref{cor:1} implies the following refinement of proposition \ref{prop:2}.
\begin{proposition}[{\bf Flocking for mono-kinetic alignment, $1<p<\nicefrac{3}{2}$}]\label{prop:3}
Let $(\rhoa,\bua)$   be a strong solution of the mono-kinetic $p$-alignment system \eqref{eqs:phydro}, $1<p<\nicefrac{3}{2}$ with   ``heavy-tailed'' communication kernel  $\phiab$ of order $0\leq \beta< \nicefrac{1}{2}$, \eqref{eq:pconnect}.
Then there is long time flocking behavior with  decay rate 
\begin{equation}\label{eq:pflocking_refind}
\begin{split}
\delE(t) \leq 
   \ds  C_RM_p (1+t)^{-\tfrac{1}{p-1}}, \qquad  1<p<\nicefrac{3}{2}, \quad 0\leq \beta< \nicefrac{1}{2}.
    \end{split}
\end{equation}
\end{proposition}
\noindent
Thus, we have $L^2_\rho$-velocity fluctuations with optimal decay rate $\lesssim (1+t)^{-\frac{1}{2p-2}}$. Moreover, integration of \eqref{eq:delV-pcontarct} with $k(D(t)) \geq C_R k(D_+)$ implies uniform decay of velocity fluctuations at the same optimal rate
\begin{equation}\label{eq:previsit}
 \max_\bx|\bua(t,\bx)-\bubar| \leq C_RM_1^{-\tfrac{2-p}{2p-2}} (1+t)^{-\tfrac{1}{2p-2}}, \qquad  1<p<\nicefrac{3}{2}.
\end{equation}
 \subsection{Agent-based description}\label{sec:agent-based}
The hydrodynamic $p$-alignment with mono-kinetic closure is the continuum counterpart of the corresponding agent-based description \eqref{eq:palignment}.
In particular, we have bounds on the velocity fluctuations  --- both the 
$\ell^2$-energy fluctuations and uniform fluctuations, which are worked out in appendix \ref{sec:pagent-based}
\begin{subequations}\label{eqs:pRic}
\begin{eqnarray}
& \qquad \ds \ddt  \delE(t)  \leq -2^{p-1}k(D(t))\left(\delE(t)\right)^p, \quad  \delE(t):=\frac{1}{2N^2}\sum_{i,j=1}^N|\bv_i(t)-\bv_j(t)|^2 \label{eq:pRica}\\
& \qquad \ds \ddt  \delbv(t)  \leq -\frac{1}{2}k(D(t))\big(\delbv(t)\big)^p, \quad \delbv(t):=\max_i|\bv_j(t)-\bvbar|, \ \bvbar(t):=\frac{1}{N}\sum_{j=1}^N \bv_j(t). \label{eq:pRicb}
\end{eqnarray}
\end{subequations}
There is one-to-one correspondence between \eqref{eqs:pRic} and   the hydrodynamic fluctuations bounds --- the $L^2_\rho$-energy fluctuations  \eqref{eq:phydroRic} and uniform velocity fluctuations in \eqref{eq:delV-pcontarct}.\newline
 When $p=1$, \eqref{eq:pRica} implies the exponential decay of heavy-tailed kernels. This should be contrasted with the case $p>1$, where  the $p$-graph Laplacian in \eqref{eq:palignment} implies polynomial  decay. 
 A typical scenario is summarized in the following proposition.

\begin{proposition}\label{pro:puniform_decay}\label{prop:4}
Consider the $p$-alignment system \eqref{eq:palignment},  with a
 ``heavy-tailed'' communication kernel  of order $0\leq\beta<1$, 
 \eqref{eq:pconnect}.
Then there is a uniform convergence towards the mean velocity
\begin{equation}\label{eq:puniform_decay}
\begin{split}
\max |\bv_i(t)-\bvbar|
  \leq \left\{ \begin{array}{ll}
  \ds C_R\, exp \big\{-C_k\myangle{t}^{(1-\beta)}\big\}\delE(0), & p=1,\\ \\
  \ds {C_R}\myangle{t}^{-\tfrac{1-\beta}{2p-2}} & p>1.
  \end{array}\right.
  \end{split}
\end{equation}
\end{proposition}

\begin{remark}[{\bf Finite-time alignment for $0\leq p <1$}]\label{rem:zalignment}
The dynamics of $p$-alignment with $p\geq1$ is driven by gradient of velocities, $\bv_j-\bv_i$. For $0\leq p <1$, the dynamics emphasizes the \emph{orientation} of velocities' gradient. The prototypical case is  $p=\nicefrac{1}{2}$, in which case \eqref{eq:palignment} reads
\begin{equation}\label{eq:zalignment}
    \left\{\begin{array}{c}
    \begin{split}
        \ddt\bx_{i}(t) &= \bv_{i}(t),   \\
        \ddt \bv_{i}(t) &= \frac{1}{\Nb}\sum_{j \neq i}^{\Nb}\phiab_{ij}(t)\frac{\bv_{j}(t)-\bv_{i}(t)}{|\bv_{j}(t)-\bv_{i}(t)|}
    \end{split}      
    \end{array}\right. \quad i=1,2,\ldots \Na.
\end{equation}
When $p=0$, \eqref{eq:zalignment}${}_2$ reads
\[
\ddt \bv_{i}(t) = \frac{1}{\Nb}\sum_{j \neq i}^{\Nb}\phiab_{ij}(t)\frac{\bv_{j}(t)-\bv_{i}(t)}{|\bv_{j}(t)-\bv_{i}(t)|^2}, \quad i=1,2,\ldots \Na.
\]
The balance of its energy fluctuations
\[
\ddt \delE(t) = -\frac{1}{2N^2} \sum_{i,j=1}^N \phi(\bx_i,\bx_j)\leq -\hf k(D(t)) \ \ \leadsto \ \ \delE(t) \leq \delE_0 - \hf\int_0^t k(D(s)){\d}s,
\]
proving that there is \emph{finite-time alignment}, $\delE(t)\stackrel{t\rightarrow t_c}{\longrightarrow}0$, for heavy-tailed kernels satisfying $\int\limits_0^t k(D(s)){\d}s \stackrel{t\rightarrow \infty}{\longrightarrow}\infty$.
Finite-time alignment is typical for  $p$-alignment in the singular range $0\leq p<1$, \cite[Theorem 2.2]{CCH2014},  
\[
\delbv(t) \leq \Big((\delbv_0)^{1-p} - 2^{p-1}(1-p)\int\limits_0^t k(D(s)){\d}s\Big)^{\tfrac{1}{1-p}}, \qquad 0\leq p <1.
\] 
 In this context, at least for $0\leq p\leq \nicefrac{1}{2}$, one encounters the need to avoid collisions, 
\[
|\bv_i(t)-\bv_j(t)|+ |\bx_i(t)-\bx_j(t)|\neq 0, \qquad i\neq j, \ t< t_c.
\]
Collision avoidance is discussed in  \cite{Mar2018} for $p\in (\nicefrac{1}{2},\nicefrac{3}{2})$ and for the  case of pure alignment, $p=1$, with possibly singulars, $k(r)=r^{-\alpha}$, in \cite{ACHL2012,Pes2014,CCH2014,CCMP2017}.
\end{remark}
\subsection{Flocking with matrix-valued communication kernel}\label{sec:matrix-kernel} 
\normalfont
Consider the alignment dynamics
\begin{subequations}\label{eqs:mhydro}
\begin{equation}\label{eq:mkernel}
          \partial_{t}(\rhoa\bua)+\nabla_\bx\cdot(\rhoa\bua\otimes\bua) =  \int \limits_{\R^d}\Phi (\bx,\by)(\bub(t,\by)-\bua(t,\bx))\rhoa(t,\bx)\rhob(t,\by)\dy, 
\end{equation}
driven by a bounded symmetric \emph{matrix} communication kernel, $\Phi (\bx,\by)=\Phi (\by,\bx)\in {\mathbb R}^{d\times d}$, of order $\beta\geq 0$
\begin{equation}\label{eq:mbeta}
    C_k \myangle{|\bx-\by|}^{-\beta}{\mathbb I}_{d\times d} \leq \Phi (\bx,\by)  \leq \phi_+{\mathbb I}_{d\times d}.
  \end{equation}
  \end{subequations}
In this case, the coupling of $\bu$-components defies a maximum principle of $\delbu(t)$ encoded in \eqref{eq:delV-contarct}. Instead, we will show below the bound $\delbu(t) \lesssim \myangle{t}^{\nicefrac{1}{2}}$. This implies   $D(t) \lesssim \myangle{t}^{\nicefrac{3}{2}}$ and hence flocking holds for heavy-tailed kernels 
 of order $\beta<\nicefrac{2}{3}$.
 To this end,  we follow our argument in the discrete setup, \cite[Proposition 3.1]{ST2021}, starting with the alignment dynamics
 \[
\partial_t\bu + \bu\cdot\nabla_\bx \bu = \int \Phi(\bx,\by)(\bu'-\bu)\rho'\dxp,
 \]
 which implies the \emph{local} energy balance
 \begin{equation}\label{eq:local-energy}
\partial_t\frac{|\bu|^2}{2} + \bu\cdot\nabla_\bx \frac{|\bu|^2}{2} = \int \big\langle \bu,\Phi(\bx,\by)(\bu'-\bu)\big\rangle\rho'\dxp.
 \end{equation}
 The integrand on the right is decomposed by polarization (suppressing time dependence)
 \[
 \begin{split}
 \langle \bu(\bx),\Phi(\bx,\by)&\big(\bu(\bx')-\bu(\bx)\big)\rangle \\
 & \equiv -\frac{1}{2}\big\langle (\bu'-\bu),\Phi(\bx,\by)(\bu'-\bu)\big\rangle
 - \frac{1}{2}\big\langle \bu,\Phi(\bx,\by)\bu\big\rangle
 +\frac{1}{2}\big\langle \bu',\Phi(\bx,\by)\bu'\big\rangle \\
 & \leq  -  C_k \myangle{|\bx-\by|}^{-\beta}\frac{|\bu|^2}{2}
 +\phi_+\frac{|\bu'|^2}{2}, \qquad \Phi(\bx,\bxp)\leq \phi_+{\mathbb I}_{d\times d}.
 \end{split}
  \]
  In the last step we used the assumed bound on $\Phi$ having a heavy-tail of order $\beta$ and satisfying a pointwise upper-bound $\phi_+$.
  Returning to \eqref{eq:local-energy} while noting that
  $\ds \int {|\bu'|^2}\rho'\dxp \leq C^2_0= 2\int \rho_0 E_0$, 
it follows that
 \begin{equation}\label{eq:boot-D}
\partial_t|\bu|^2 + \bu\cdot\nabla_\bx |\bu|^2\leq  -C_k\myangle{D(t)}^{-\beta}M|\bu|^2+\phi_+C^2_0.
 \end{equation}
 By the maximum principle (--- we ignore the dissipative term on the right),
 \[
 |\bu(t,\cdot)|^2\leq \max|\bu_0|^2+C' t, \qquad C':=\phi_+C^2_0,
 \]
  and hence \eqref{eq:sizeofDt} holds with $\gamma=\nicefrac{3}{2}$, in view of 
\begin{equation}\label{eq:D-expandb}
 \ddt {D}(t)\leq 2\max|\bu(t,\cdot)| \ \ \leadsto \ \ D(t) \leq D_0 +\frac{4}{3C'}\big(\max|\bu_0|^2+C't\big)^{\nicefrac{3}{2}}.
 \end{equation} 
 We can now use a bootstrap argument: starting with $\gamma=\nicefrac{3}{2}$ we insert the bound  $D(t)\lesssim \myangle{t}^\gamma$ of \eqref{eq:D-expandb} into the right side of \eqref{eq:boot-D} and we have the maximum bound
 \[
 |\bu(t,\cdot)|^2\leq \max|\bu_0|^2+C' (1+t)^{\beta\gamma} \ \ \leadsto \ \ D(t) \leq D_0 +\frac{2}{C'\gamma'}\big(\max|\bu_0|^2+C'(1+t)\big)^{\gamma'}, \quad \gamma'=1+\nicefrac{\beta\gamma}{2}.
 \]
 Iterating,  $\gamma\mapsto \gamma'$, we end up with a fixed point $\gamma=\frac{2}{2-\beta}$, and with the improved bounds, still in the range of $\beta<\nicefrac{2}{3}$,
 \[
 |\bu(t)| \lesssim \myangle{t}^{\tfrac{\beta}{2-\beta}}, \qquad D(t) \leq C_D\myangle{t}^{\tfrac{2}{2-\beta}}, \quad \beta<\nicefrac{2}{3}.
 \]
 Corollary \ref{cor:1} implies the following.
 \begin{proposition}[{\bf Flocking for matrix-based alignment}]\label{prop:5}
Let $(\rhoa,\bua)$   be a strong solution of the hydrodynamic  alignment  \eqref{eqs:mhydro} with
 ``heavy-tailed''  matrix communication kernel  $\Phi$ of order $ \beta< \nicefrac{2}{3}$.
There is long time flocking behavior with  fractional exponential decay rate 
\begin{equation}\label{eq:mflocking}
\delE(t)   \leq 
  C_R\, exp\big\{- M_1 t^{\tfrac{2-3\beta}{2-\beta}}\big\}\delE(0).
\end{equation}
\end{proposition}

\section{Flocking of hydrodynamic $p$-alignment with entropic pressure}\label{sec:with-pressure}
We consider hydrodynamic alignment \eqref{eqs:phydro} driven by the class of \emph{singular kernels}
$k_p(r):=r^{-(d+2sp)}, 0<s<1, \ p\geq 1$,
\[
\partial_t(\rho \bu)+\nabla_\bx(\rho \bu \otimes \bu+\Pressure) = 
p.v. \int \limits_{\Omegaa(t)} \frac{|\bu'-\bu|^{2p-2}\big(\bu(t,\bx')-\bu(t,\bx)\big)}{|\bx'-\bx|^{d+2sp}}\rho(t,\bx)\rho(t,\bxp)\dxp, \quad 0<s<1.
\]
We emphasize that in this case of strongly singular kernels, there is no formal justification for the passage from the  agent-based description \eqref{eq:palignment} to the hydrodynamic description. In particular, the near-origin integrability sought in \eqref{eq:pconnect} is given up for the usual notion of singular integration in terms of principle value ($p.v.$). The alignment term on the right  amounts to a \emph{weighted} fractional $2p$-Laplacian, $(-\Delta)^s_{2p}$, which is  properly interpreted to act  on $\textnormal{supp} \, \rho(t,\cdot)$; see \cite{TGCV2021,BV2015} and the references therein. 
 
The tail of the singular kernel, $k_p(r)=r^{-(d+2sp)}, r\gg R$, is too thin to enforce the heavy tail condition sought in corollary  \ref{cor:1}. Accordingly, we keep the singular `head' and adjust it with  the ``heavy tail'' or order $\beta$
\begin{equation}\label{eq:multiDsingdyn}
\phiab_{s,\beta} (\bx,\bx')\left\{\begin{array}{ll} = |\bx-\bx'|^{-(d+2sp)}, & |\bx-\bx'|\leq R \ \ \mbox{with} \ \ 0< s< 1\\ \\
 \geq C_k\myangle{|\bx-\bx'|}^{-\beta}, & |\bx-\bx'|> R
 \end{array}\right.
\end{equation}
Clearly, there exists  a constant, $K=K_R$, such that  $k_p(|\bx-\bxp|) \leq K_R\phi_{s,\beta}(\bx,\bxp)$ for all $(\bx,\bxp)$. Without loss of generality, we  may assume that the spatial scale  $R$  is large enough, 
$\ds  {(1+R)^\beta}R^{-(d+2sp)} <C_k$, so that we may take $K_R=1$,
\begin{equation}\label{eq:phi-dominates}
k_p(|\bx-\bxp|) \leq \phi_{s,\beta}(\bx,\bxp), \qquad \forall \bx,\bxp \in \R^d.
\end{equation} 
We refer to such heavy-tailed, singular kernels as having order $(s,\beta)$. 
If we let $\phi_\beta$ denote its tail of order $\beta$ then the  $p$-alignment dynamics now  reads
\begin{equation}\label{eq:sing-palignment}
\begin{split}
\partial_t(\rho \bu)&+\nabla_\bx\cdot(\rho\bu\otimes\bu +\Pressure) \\
   & =  p.v. \!\!\!\int \limits_{|\by-\bx|\leq R} \hspace*{-0.3cm} \frac{|\bu'-\bu|^{2p-2}\big(\bu'-\bu'\big)}{|\by-\bx|^{d+2sp}}\rho\rho'\dy \\
 &  +\hspace*{-0.4cm}\int \limits_{|\by-\bx|> R} \hspace*{-0.4cm}\phi_{\beta}(\bx,\bxp)|\bu'-\bu|^{2p-2}\big(\bu'-\bu'\big)\rho\rho'\dy, \quad \phi_\beta(\bx,\bxp)\geq C_k(1+|\bx-\bxp|)^{-\beta}.
 \end{split}
\end{equation}
\begin{remark}[{\bf Entropic pressure with singular kernel}]
In case of singular kernel $\phi_{s,\beta}$, we need to adjust the definition \ref{def:meso-ppressure} of entropic pressure, 
 \begin{equation}\label{eq:sing-ppressure}
 \partial_{t}(\rhoa \ea )+\nabla_{\bx}\cdot(\rhoa \ea \bua +\bqa )+\textnormal{trace}(\Pressurea\nabla\bua) \leq - \hf k_p(D(t))\int \limits_{{\mathcal S}(t)} \big((2\ea)^p+(2\ebp)^p\big) \rho\rhobp \dy.
 \end{equation}
  Thus, the entropic part of the internal energy  avoids the singularity of $\phi_{s,\beta}$ and emphasizes only its tail behavior. It leads to the `adjusted' energy fluctuations bound 
 \begin{equation}\label{eq:energyvssing}
\ddt \delE(t)
    \leq 
- \hf \iint \limits_{{\mathcal S}(t)\times {\mathcal S}(t)}\left\{\phi_{s,\beta} (\bx,\by) |\bua -\bubp |^{2p}+k_p(D(t))\big((2\ea)^p+(2\eb)^p\big)\right\}\rhoa\rhobp \dx\dy,
 \end{equation}
 which in turn, arguing along the  lines of \eqref{eq:delEpdecay}, yields \eqref{eq:phydroRic}; that is, the main theorem \ref{thm:main1} and its corollary \ref{cor:1} survive.
 In particular,  the enstrophy bound \eqref{eq:ens} holds \myr{for $\phi=\phi_{s,\beta}$. Taking into account  \eqref{eq:phi-dominates}, $\phi_{s,\beta}(\bx,\bxp)\geq |\bxp-\bx|^{-(d+2sp)}$, we find}
 \begin{equation}\label{eq:end-revisited}
\myr{\ddt \delE(t)
    \leq 
- \hf \iint \limits_{{\mathcal S}(t)\times {\mathcal S}(t)} \frac{|\bua(t,\bxp) -\bu(t,\bx)|^{2p}}{|\bxp-\bx|^{d+2sp}}\rhoa\rhobp \dx\dy.}
 \end{equation}
 
\end{remark}
The presence of pressure, let alone a pressure with an `unknown' closure, couples the different components of velocity in a manner that defies a straightforward derivation of a uniform bound on velocity fluctuations, $\delbu(t)$, along the lines of what we have done in the mono-kinetic case.  Instead, 
we introduce a new strategy for verifying flocking in this case, in which we use an enstrophy bound associated with the  singular kernel, 
$k_p(r)=r^{-(d+2sp)}$, in order to  control  the diameter $D(t)\lesssim \myangle{t}^{\gamma}$. This enables us to treat the flocking in presence of entropic pressure.
The remarkable aspect here is that although the presence of pressure   defies  a maximum principle on the velocity field,  the corresponding   enstrophy bound associated with \eqref{eq:sing-palignment} will suffice for control of velocity fluctuations and  hence flocking will follow. Thus, short-term interactions governed by kernel with a \emph{singular head}  secure the spread of velocity fluctuations, while \emph{heavy-tailed}  kernel governing the long-term  interactions  secure flocking. 

\subsection{Enstrophy and dispersion bounds}\label{sec:ens-and-dispersion}
 Throughout this section we make the following assumptions.
 
 \begin{enumerate}[itemindent=.01cm,label={(H\textit{\arabic*})}]
 \item\label{H1} The alignment hydrodynamics \eqref{eq:hydro},\eqref{eq:sing-palignment} admits a strong entropic solution, \eqref{eq:sing-ppressure}. 
\item\label{H2} The support,   $\Omegaa(0)=\textnormal{supp}\,\rhoa(0,\cdot)$, has a smooth boundary satisfying a Lipschitz or a cone condition. 

\item\label{H3} The dynamics remains uniformly bounded away from vacuum, namely --- there exists $\rho_->0$ such that  
\[
\min_{\bx\in \Omegaa(t)}\rhoa(t,\bx)\geq \rho_->0, \quad t\geq 0.
\]

\end{enumerate}
Several comments regarding these assumptions are in order. The literature about the question of global regularity, \ref{H1},   is devoted mostly to mono-kinetic ``pressure-less'' closure; we mention the one-dimensional studies \cite{TT2014,CCTT2016,HT2017,ST2017a, ST2017b,ST2020b,Tan2021,LS2022}, the two-dimensional case \cite{HT2017} and multi-dimensional  cases \cite{Shv2019, DMPW2019, CTT2021,Tad2022b}. Much less is known about alignment with pressure, typically when (scalar) pressure is augmented with additional process of relaxation and/or dissipation, \cite{Cho2019,CDS2020,TCGW2020}. On the other hand, there are relatively few works on weak solutions of \eqref{eqs:hydro}, \cite{CCR2011,CFGS2017,LT2021}. As for \ref{H2}, we are aware of only few results on the geometric structures that emerge from alignment, \cite{LS2019,LLST2022}. The question of uniform bound away from vacuum  assume in \ref{H3} plays an important role in driving global regularity \cite{Tan2020,Shv2021,AC2021a,Tad2021}. It can be relaxed to allow mild time decay, e.g., $\rho_-(t) \gtrsim \myangle{t}^{-\nicefrac{1}{2}}$, \cite[Theorem 1.1]{ST2020b}, \cite[Theorem 3]{Tad2021}, but as already noted in previous works, some sort of non-vacuous assumption is necessary.  \newline 

 We begin by noting that since  $\phi_{s,\beta}$ dominates $k_p(r)$, \eqref{eq:phi-dominates}, then by the non-vacuous hypothesis \ref{H3}, $\rho\geq \rho_->0$, we   have the Sobolev bound 
 \[
 \begin{split}
 \iint \limits_{\Omegaa(t)\times \Omegaa(t)}\frac{|\bu(t,\by)-\bu(t,\bx)|^{2p}}{|\by-\bx|^{d+2sp}}\dx\dy  \leq  C_\rho^2\iint \limits_{\Omegaa(t)\times \Omegaa(t)}\frac{|\bu(t,\by)-\bu(t,\bx)|^{2p}}{|\by-\bx|^{d+2sp}}\rho\rho'\dx\dy, \quad C_\rho:=\frac{1}{\rho_-}.
   \end{split}
   \]
The space-time  enstrophy  bound \eqref{eq:ens}, or more precisely --- its singular version in \eqref{eq:end-revisited}, then yields
\begin{equation}\label{eq:mDsingHs}
 \int \limits_0^t \|\bua(\tau,\cdot)\|^{2p}_{{}_{\dot{W}^{s,2p}(\Omegaa)}}\d\tau\leq C^2_\rho  C^2_0, \quad  \|\bu(t,\cdot)\|^{2p}_{{}_{\dot{W}^{s,2p}(\Omegaa)}}:=\iint \limits_{\Omegaa(t)\times \Omegaa(t)}\frac{|\bua(t,\by)-\bua(t,\bx)|^{2p}}{|\by-\bx|^{d+2sp}}\dx\dy.
 \end{equation}
  The enstrophy bound \eqref{eq:mDsingHs}   guarantees that the  velocity $\bua$ slows down the dispersion of the crowd so that its diameter $D(t)$ may not grow faster than $\lesssim \myangle{t}^{\gamma}$. Below we derive sharp bounds on the dispersion rate $\gamma$.\newline
  To this end, we note that propagation along particles paths in \eqref{eq:hydro}${}_1$ yields, as in \eqref{eq:mono-dispersion},
  \[
  \ddt D(t) \leq \delbu(t), \quad \delbu(t)=\max_{\bx,\bxp\in {\mathcal S}(t)}|\bu(t,\bxp)-\bu(t,\bx)|.
  \]
 By Gagliardo-Nirenberg inequality (which we recall in appendix \ref{sec:GN} below),
\begin{equation}\label{eq:GNS}
 |\bu(t,\bx)-\bu(t,\bx')|\leq C_s\|\bu\|_{{}_{\dot{W}^{s,2p}({\mathcal S}(t))}}|\bx-\bx'|^{s-\theta}, \qquad \bx,\bx'\in \Omegaa(t), \ \ \theta:=\frac{d}{2p} < s<1.
\end{equation}
This yields, 
$\ds \ddt \Da(t)\leq \delbu(t) |\leq C_s\|\bu(t,\cdot)\|_{\dot{W}^{s,2p}({\mathcal S}(t))}D^{s-\theta}(t)$, or
\begin{equation}\label{eq:dynamicD}
\ddt \Da^{1+\theta-s}(t) \leq C'_s\|\bu(t,\cdot)\|_{\dot{W}^{s,2p}({\mathcal S}(t))}, \qquad C'_s=\myr{(1+\theta-s)}C_s,
\end{equation}
and hence, in view of \eqref{eq:mDsingHs},
\[
\begin{split}
\Da^{1+\theta-s}(t)  & \leq \Da^{1+\theta-s}_0  \\
 & \ \ \ + \Big(\int \limits_0^t \|\bu(\tau,\cdot)\|^{2p}_{\dot{W}^{s,2p}({\mathcal S}(t))}{\d}\tau\Big)^{\tfrac{1}{2p}}\Big(\int \limits_0^t 1{\d}\tau\Big)^{\tfrac{1}{(2p)'}}  \leq \Da^{1+\theta-s}_0 
+ C'_s( C_\rho  C_0)^{\tfrac{1}{p}}\, t^{\tfrac{1}{(2p)'}}.
\end{split}
\]
We conclude that the crowd of multi-dimensional $p$-alignment dynamics, 
\eqref{eq:sing-palignment} can be dispersed at a rate no faster than
\begin{equation}\label{eq:Dabeta}
\Da(t) \leq C_D(1+t)^{\gamma_p}, \qquad \gamma_p= \frac{2p-1}{2p(1+\theta-s)}, \quad \theta=\frac{d}{2p}<s <1.
\end{equation}
This bound can be improved \myr{using a bootstrap argument outlined in appendix \ref{sec:dispersion}}. In particular, for $1<p<\nicefrac{3}{2}$ we  obtain  a \emph{uniform} dispersion bound which we summarize in the following key result.
\begin{lemma}[{\bf Uniform dispersion bound  for $p$-alignment with singular kernels}]\label{lem:dispersion}
Consider the multi-dimensional $p$-alignment dynamics, 
\eqref{eq:sing-palignment}, $1<p<\nicefrac{3}{2}$, with heavy-tailed, singular kernel of order $(s,\beta)$, satisfying \ref{H1}--\ref{H3}. 
Then we have a uniform bound 
\begin{equation}\label{eq:Duniform}
\Da(t) \leq D_+, \qquad 0\leq  \myr{\beta< (\nicefrac{3}{2}-p)d}, \quad 1<p<\nicefrac{3}{2}. 
\end{equation}
\end{lemma}
\begin{remark}\label{rem:more-on-p}
Observe that since we require $\myr{d=2p\theta<3}$, the uniform bound \eqref{eq:Duniform} is restricted to one- and two-dimensional cases.\newline
We are unable to secure such a uniform dispersion bound for $p>\nicefrac{3}{2}$, but we can still improve the dispersion bound  \eqref{eq:Dabeta} as shown in  remark \ref{rem:pgt32} below,
\[
D(t)\leq C'_D(1+t)^\gamma, \qquad  \gamma= \frac{2p\big(p-\nicefrac{3}{2}\big)}{(p-1)d-\beta},   \qquad 0\leq \beta < \frac{d}{2p-1}, \quad  p>\nicefrac{3}{2}.
\]
\end{remark}
\subsection{Flocking of alignment with pressure. The one-dimensional case}
The case of pure alignment $p=1$, restricts the use of lemma \ref{lem:dispersion}  to the one-dimensional case ($d<2\myr{p}$),
 driven by singular kernel   $k_1(r)=r^{-(1+2s)}, \hf<s<1$, with  $\beta$-tailed adjustment
\begin{subequations}\label{eqs:1Dsingdyn}
\begin{equation}\label{eq:1Dsingdyn}
\begin{split}
\partial_t(\rho u)+\partial_x(\rho u^2+\pressure) & = p.v. \hspace*{-0.5cm}\int \limits_{|x'-x|\leq R}\hspace*{-0.3cm}\frac{u(t,x')-u(t,x)}{|x-x'|^{1+2s}}\rho(t,x'){\d}x' \\
 & \ \ \ + \hspace*{-0.5cm} \int \limits_ {|x'-x|> R}\hspace*{-0.3cm}\phi_{\beta}(x,x')\big(u(t,x')-u(t,x)\big)\rho(t,x'){\d}x'.
\end{split}
\end{equation}
The integrals on the right are restricted to the interval  ${\mathcal S}(t)=[\rho_-(t),\rho_+(t)]$ supporting $\rho(t,\cdot)$, $\phi_{\beta}$ is a $\beta$-tailed communication kernel,
\begin{equation}
\phi_{\beta}(x,\xp) \geq C_k(1+|x-x'|)^{-\beta}, \qquad |x-x'|\geq R,
\end{equation}
and $\pressure$ is any scalar entropic pressure satisfying \eqref{eq:meso-pressure}, or more precisely --- its singular version \eqref{eq:sing-ppressure},
\begin{equation}\label{eq:1Dsingpres}
\partial_t (\rho \pressure) +\partial_x(\rho \pressure u+q) +2\pressure \partial_x u \leq -2\pressure D^{-(1+2s)}(t)M.
\end{equation}
\end{subequations}
By  \ref{eq:Dabeta} we can apply corollary \eqref{cor:1} with $\gamma_1=\tfrac{1}{3-2s}$ which  yields the following.
 \begin{theorem}[{\bf One-dimensional alignment, $p=1$}]\label{thm:main3}
Consider the one-dimensional alignment dynamics 
\eqref{eqs:1Dsingdyn} and assume  \ref{H1},\ref{H3}, hold. Let $(\rho,u,\pressure)$   be a strong entrpoic solution
 with $\beta$-tailed singular kernel, $\phi_{\beta}$, satisfying the heavy-tail condition 
 \begin{equation}\label{eq:heavys}
 \beta+2s <3,\qquad \beta\geq 0, \quad \hf<s<1.
  \end{equation}
   Then there is a large time flocking behavior  with fractional exponential rate
\begin{equation}\label{eq:1Dflocking_decay}
\delE(t) 
  \leq C_R \,exp\Big\{-2MC_k\myangle{t}^{\tfrac{3-2s-\beta}{3-2s}}\Big\}\delE(0).
\end{equation}
\end{theorem}
\noindent
This extends the  mono-kinetic `pressure-less' studies in \cite{ST2017a, ST2017b, ST2018a, DKRT2018, DMPW2019, MMPZ2019}.
 It is instructive to compare this result with flocking statement in  the mono-kinetic closure, which is based on the uniform bound on velocity, $D(t) \lesssim \myangle{t}$. Theorem \ref{thm:main3} allows for a \emph{larger} class of heavy-tailed kernels since  it is based on a sharper bound on the velocity fluctuations, leading to $D(t) \lesssim \myangle{t}^{\gamma}$ with $\gamma <1$. \myr{This result can be further improved by extending the uniform dispersion bound in lemma \ref{lem:dispersion} to the limiting case $p=1$.}
 
 \subsection{Flocking of $p$-alignment with pressure. The multi-dimensional case}

 We consider the $p$-alignment  dynamics \eqref{eq:sing-palignment} driven by singular kernel   $k_p(r)=r^{-(d+2sp)}, \tfrac{d}{2p}<s<1$.
 Using \eqref{eq:Dabeta} we can apply corollary \ref{cor:1} with $\gamma=\gamma_p$ which  yields the following.
 \begin{theorem}[{\bf Multi-dimensional alignment, $p>1$}]\label{thm:main4}
Consider the multi-dimensional $p$-alignment dynamics 
\eqref{eq:sing-palignment} and assume \ref{H1}--\ref{H3} hold. Let $(\rho,\bu,\Pressure)$   be a strong entrpoic solution, \eqref{eq:sing-ppressure},  
  with a $\beta$-tailed singular kernel $\phi_{s,\beta}$, satisfying the heavy-tail condition
  \[
  \beta\gamma_p<1, \qquad \beta\geq0, \quad \gamma_p:=\frac{2p-1}{2p(1+\theta-s)}, \quad \theta=\frac{d}{2p}<s<1.
  \] 
Then there is a large time flocking behavior  with polynomial decay rate of order
\begin{equation}\label{eq:mDflocking_decay}
\delE(t) 
  \leq C_RM_p\,t^{-\tfrac{1-\beta\gamma_p}{p-1}}.
\end{equation}
\end{theorem}
\begin{remark}[{\bf Decay of internal fluctuations}]
A sufficient condition for the heavy-tailed restriction $\beta\gamma_p<1$ sought in \eqref{eq:mDflocking_decay} is given by 
 \begin{equation}\label{eq:mDheavys}
 \beta \leq \frac{d}{2p-1} \ \ \leadsto \ \  \myr{\beta\gamma_p<\beta\frac{2p-1}{d}\leq 1}.
 \end{equation}
It still allows heavy-tails of order $\beta\geq1$,
compared with the $\beta<1$ restriction in the mono-kinetic closure.
In particular, when $\beta=\frac{d}{2p-1}$ one finds the decay  of order
\[
\big(\delE(t)\big)^{p-1} \lesssim t^{-(1-\beta\gamma_p)} \lesssim t^{-  \tfrac{1-s}{1+\theta-s}}.
\]
\end{remark}
\begin{remark}
Theorem \ref{thm:main4} implies the decay of both --- the macroscopic velocity fluctuations
$\ds \int |\bu-\bubar|^2\rho\dx$ and, in the context of kinetic formulation, the microscopic fluctuations $\ds \iint |\bv-\bu|^2 f_N \dv\dx$.
\end{remark}
 The decay bound \eqref{eq:mDflocking_decay} is not sharp, a reflection of the fact that the dispersion bound \eqref{eq:Dabeta} can be improved with smaller $\gamma_p$ (\myr{as noted in remark  \ref{rem:more-on-p}}). In particular, when $p$ is in the restricted range $1<p<\nicefrac{3}{2}$, then   corollary \ref{cor:1} applies  with $\gamma=0$ and $C_D=D_+$ which  yields the following.
\begin{theorem}[{\bf Multi-dimensional alignment $1<p<\nicefrac{3}{2}$}]\label{thm:main5}
Consider the multi-dimensional $p$-alignment dynamics 
\eqref{eq:sing-palignment}, $1<p<\nicefrac{3}{2}$ and assume \ref{H1}--\ref{H3} hold. Let $(\rho,\bu,\Pressure)$   be a strong entrpoic solution, \eqref{eq:sing-ppressure},  
  with a heavy-tailed singular kernel of order $(s,\beta)$. 
Then there is a large time flocking behavior  with polynomial decay rate of order
\begin{equation}\label{eq:pDflocking_refind}
\begin{split}
\delE(t) \leq 
   \ds  C_RM_p (1+t)^{-\tfrac{1}{p-1}}, \qquad   0\leq \beta< \myr{(\nicefrac{3}{2}-p)}d, \quad \frac{d}{2p}<s<1, \quad 1<p<\nicefrac{3}{2}.
    \end{split}
\end{equation}
\end{theorem}
\noindent
Theorem \ref{thm:main5} is the analogue of the mono-kinetic ``pressure-less'' case in proposition \ref{prop:3}. In particular, it is rather remarkable that we obtain the same optimal decay rate of order $\tfrac{1}{p-1}$ in the respective range $1<p<\nicefrac{3}{2}$ for the one- and two-dimensional cases. An optimal flocking scenario with a uniform dispersion bound  remains open for $d\geq 3$.

\appendix
\section{Derivation of entropic inequality in $p$-alignment}\label{sec:derivation}
\subsection{From agent-based to hydrodynamic description}\label{sec:hydro-des} 
We begin with the  passage from  the agent-based dynamics of $p$-alignment \eqref{eq:palignment} to its hydrodynamic description \eqref{eqs:phydro}.
The large crowd dynamics is encoded in terms of their empirical distribution $\displaystyle \fa(t,\bx,\bv):= \frac{1}{\Na}\sum \limits_{i=1}^{\Na} \delta_{\bx_i(t)}(\bx)\otimes \delta_{\bv_i(t)}(\bv)$, which are governed by the kinetic Valsov equation in state variables   $(t,\bx,\bv)\in \R_t\times\R^d\times \R^d$, 
\begin{equation}\label{eq:Q}
\partial_t \fa +\bv\cdot\nabla_\bx \fa + \nabla_\bv\cdot \Qp(\fa,\fa)=0,
\end{equation}
and driven by interaction kernel
\[
\Qp(\fa,\fb)\myr{(t,\bx,\bv)} := \int\limits_{{\mathcal S}(t)\,\,}\myr{\int \limits_{\R^d}}\phiab (\bx,\bxp)|\bvp-\bv|^{2p-2}(\bvp-\bv)\fa \fbp \dvp\dxp.
\]
We distinguish between the cases   of `pure' alignment, $Q_1=Q$, and enhanced $p$-alignment $Q_p$ of order $p>1$.\newline 
For $p=1$, the large crowd dynamics of  $\fa$'s is captured by their first two moments which we assume to exist --- the density $\ds \rho:=\lim_{\Na\rightarrow \infty} \int_{\R^d} \fa(t,\bx,\bv)\dv$ and momentum $\ds \rho \bu:=\lim_{\Na\rightarrow \infty} \int_{\R^d} \bv \fa(t,\bx,\bv)\dv$; that is,  
\begin{equation}\label{eq:dandm}
 \rhoa(\bvp-\bu)=\lim_{\Na\rightarrow \infty} \int \limits_{\R^d} (\bvp-\bv)\fa(t,\bx,\bv)\dv \ \ \textnormal{for all} \ \ \bvp\in \R^d.
 \end{equation}
 Integration of \eqref{eq:Q} yields the mass equation \eqref{eq:hydro}${}_1$ 
\begin{subequations}\label{eqs:hydro-revisited}
\begin{equation}
\partial_{t}\rhoa+\nabla_\bx\cdot(\rhoa\bua) = 0.
\end{equation}
The first $\bv$-moment of \eqref{eq:Q} yields 
\[
\partial_t \int \limits_{\R^d}\bv \fa\dv  =-\nabla_\bx\cdot \int \limits_{\R^d}\bv \otimes \bv \fa\dv  +
  \int \limits_{\R^d} Q_{\myr{1}}(\fa,\fa)\dv.
  \]
  We now treat the two terms on the right. For the first term, we decompose
  $\bv \otimes \bv \equiv  -\bua\otimes \bua+(\bv\otimes\bu+\bua\otimes\bv)+ (\bv-\bua)\otimes (\bv-\bu)$, where the corresponding first two moments of $\fa$ add up to  $\bua \otimes(\rho\bua)= \rhoa\bua \otimes \bua$, while the third yields the pressure tensor \eqref{eq:pressure}, 
 \[
 \lim_{\Na\rightarrow \infty}\int \limits_{\R^d} \bv \otimes\bv \fa\dv =
 \rhoa\bua \otimes \bua +  \Pressurea, \qquad \Pressurea=\lim_{\Na\rightarrow \infty}\int \limits_{\R^d}(\bv-\bua)\otimes(\bv-\bua) \fa.
 \]
 The second term on the right yields
 \[
 \lim_{\Na\rightarrow \infty}\int \limits_{\R^d}Q_1(\fa,\fa)\dv = \int  \limits_{{\mathcal S}(t)}\phiab (\bx,\bxp)\left(\rhobp\bubp\rhoa-\rhoa\bua\rhobp \right)\dxp =\aligna(\rho,\bu),
\]
and we recover the momentum equation \eqref{eq:hydro}${}_2$
\[
 \partial_{t}(\rhoa\bua)+\nabla_\bx\cdot(\rhoa\bua\otimes\bua+\Pressurea) = \aligna(\rho,\bu).
\]
For $p>1$\myr{,} we \emph{assume} the existence of the corresponding higher moments (which are compatible with the mono-kinetic Maxwellian \eqref{eq:mono-closure}),
\[
\rhoa |\bvp-\bua|^{2p-2}(\bvp-\bua):=\lim_{\Na\rightarrow \infty} \int \limits_{\R^d} |\bvp-\bv|^{2p-2}(\bvp-\bv) \fa(t,\bx,\bv)\dv,
\]
in which case  the interaction kernel yields,
\[
\lim_{\Na\rightarrow \infty}\int \limits_{\R^d}\Qp(\fa,\fa)\dv = \int  \limits_{{\mathcal S}(t)}\phiab (\bx,\bxp)|\bubp-\bu|^{2p-2}\left(\rhobp\bubp\rhoa-\rhoa\bua\rhobp \right)\dxp =\aligna_p(\rho,\bu),
\]
and we recover the momentum equation \eqref{eq:phydro}${}_2$
\begin{equation}
 \partial_{t}(\rhoa\bua)+\nabla_\bx\cdot(\rhoa\bua\otimes\bua+\Pressurea) = \aligna_p(\rho,\bu).
\end{equation}
\end{subequations}
\myr{
In fact, we are not restricted here by the mono-kinetic closure assumption:
 for any kinetic closure we have 
\[
\begin{split}
\int \limits_{{\mathcal S}(t)} \int \limits_{\R^d} Q_p(\fa,\fa)(t,\bx,\bv)&\dv\dx
\\
 & =\iint \limits_{{\mathcal S}(t)\times {\mathcal S}(t)} \ \,\iint \limits_{\R^d\times \R^d} \phi(\bx,\bxp) |\bvp-\bv|^{2p-2} (\bvp-\bv)\fa'\fa\dvp\dv\dx'\dx=0.
 \end{split}
\]
This follows by the anti-symmetry of the integrand on the right, and hence the zero-average condition for $p$-alignment sought in remark \ref{rem:zero-average-Ap} holds,
\begin{equation}\label{eq:zero-kaverage-Ap}
\int \limits_{{\mathcal S}(t)} \aligna_p(\rho,\bu)(t,\bx)\dx=
\lim_{\Na\rightarrow \infty}\int \limits_{{\mathcal S}(t)}  \int \limits_{\R^d}\Qp(\fa,\fa)(t,\bx,\bv)\dv\dx=0.
\end{equation}
}
Observe that system \eqref{eqs:hydro-revisited} is not a purely hydrodynamic description, since the  pressure,   $\Pressure$,   still requires a  closure of the  second-order moments of $\fa$. Thus,   the    alignment dynamics in \eqref{eqs:hydro-revisited} is  left open at the mesoscale, subject to the notion of entropic pressure in definition \ref{def:meso-pressure} for $p=1$ and definition \eqref{def:meso-ppressure} for $p > 1$.

\subsection{Entropic pressure in kinetic formulation of $p$-alignment}\label{sec:meso-ppressure}
We follow  the balance of the internal energy balance  as preparation for studying the  large-time behavior of `pure'  hydrodynamic alignment, $p=1$, in \eqref{eqs:hydro} and hydrodynamic $p$-alignment, $p>1$, in \eqref{eqs:phydro}. The total energy is given by the second moment which is assumed to exist
\begin{equation*}
    \rhoa \Ea(t,\bx) = \lim_{\Na  \rightarrow\infty}\int \limits_{\R^d}\frac{1}{2}|\bv|^{2}\fa(t,\bx,\bv)\dv.
\end{equation*}
It is decomposed into kinetic and internal energy corresponding to the decomposition $\displaystyle \hf|\bv|^{2} \equiv \hf|\bua|^{2} + \hf|\bv-\bua|^{2}+\bua\cdot(\bv-\bua)$. Noting that  $\ds \int_{\R^d}(\bv-\bua)\fa\dv=0$,  we find
\begin{equation*}
\rhoa \Ea = \frac{\rhoa}{2}|\bua|^{2}+\rhoa \ea, \qquad \rhoa \ea:= \lim_{\Na  \rightarrow\infty}\hf\int \limits_{\R^d} |\bv-\bua|^2\fa\dv.
\end{equation*}
The balance of internal energy, $\rhoa\ea$,  is obtained by integrating (\ref{eq:Q}) against $\displaystyle \frac{|\bv-\bu|^{2}}{2}$, which yields 
\begin{equation}\label{eq:kinternala}
\partial_t (\rhoa \ea) + \int \limits_{\R^d} \frac{|\bv-\bu|^{2}}{2}\bv\cdot\nabla_{\bx}\fa\dv = \int \limits_{\R^d} (\bv-\bu)\cdot \Qp(\fa, \fb)\dv.
\end{equation}
The integral  on the left can be expressed as a perfect divergence of the cubic moments $\ds \bqa_N:=\hf\int |\bv-\bua |^{2}(\bv-\bua )\fa\dv$ (all integrals are taken over $\R^d$)
\[
\begin{split}
  \int   \frac{|\bv-\bu|^{2}}{2}&\bv\cdot\nabla_{\bx}\fa\dv \\
  & =  \nabla_{\bx}\cdot \int\frac{|\bv-\bu|^{2}}{2}\bv \fa\dv
   - \int \bv\cdot\nabla_\bx \frac{|\bv-\bu|^{2}}{2} \fa\dv \\
  &  = \nabla_{\bx}\cdot \int\frac{|\bv-\bu|^{2}}{2}\big(\bua +(\bv-\bua)\big) \fa\dv
  + \sum_{i,j}\int v_j(v_i-u_i)\frac{\partial u_i}{\partial x_j}\fa \dv  \\
   & =  \nabla_{\bx}\cdot \Big(\int\frac{|\bv-\bu|^{2}}{2}\fa\dv\, \bua +\bqa_N\Big)
  + \sum_{i,j}\int (v_j-u_j)(v_i-u_i)\fa\dv\frac{\partial u_i}{\partial x_j}.
  \end{split}
  \]
    Taking the limit we find the term
$\nabla_{\bx}\cdot (\rho\ea\bu +\bqa_h) +\textnormal{trace}(\Pressurea\nabla\bu)$, 
with heat-flux,
\begin{equation}\label{eq:heat-flux}
\ds \bqa_h :=\lim_{\Na\rightarrow \infty} \hf\int \limits_{\R^d} |\bv-\bua |^{2}(\bv-\bua )\fa\dv,
\end{equation}
 and \eqref{eq:kinternala} yields
\begin{equation}\label{eq:kinternalb}
\partial_t (\rhoa \ea) + \nabla_{\bx}\cdot (\rho\ea\bu +\bqa_h) +\textnormal{trace}(\Pressurea\nabla\bu) = \lim_{\Na\rightarrow \infty}\int \limits_{\R^d} (\bv-\bu)\cdot \Qp(\fa, \fb)\dv.
\end{equation}
It remains to consider the  moment of the alignment-based  term  on the right. We distinguish between the cases $p=1$ and $p>1$.\newline
{\bf The case  $p=1$}. We  split $\bv-\bvp\equiv(\bv-\bu)+(\bu-\bvp)$,
 \begin{equation}\label{eq:casepeq1}
 \begin{split}
\int \limits_{\R^d} (\bv-\bu)\cdot \Qa(\fa, \fb)\dv  & = -\int \limits_{{\mathcal S}(t)}\phiab (\bx,\bx')\iint \limits_{\R^d\times \R^d}(\bv-\bu)\cdot(\bv-\bvp)\fa \fbp \dvp\dv\dy \\
   & =   -\int \limits_{{\mathcal S}(t)}\phiab (\bx,\bx')\iint \limits_{\R^d\times \R^d}|\bv-\bu|^2\fa \fbp\dvp\dv\dy \\
   & \ \ \  - \int \limits_{{\mathcal S}(t)}\phiab (\bx,\bx')\int \limits_{\R^d}(\bv-\bu)\fa \dv \cdot \int \limits_{\R^d}(\bu-\bvp)\fbp\dvp\dy.  
\end{split}
\end{equation}
 The the first  integral on the right ends up with 
    $\ds -2\int \limits_{\Omegab(t)} \phiab(\bx,\by)\ea\rho\rhobp\dy$, and since the second integral on the right vanishes, \eqref{eq:kinternalb} now reads  \cite{HT2008}\footnote{This corrects a series of typos in our statement of \cite[Lemma 5.1]{HT2008}} 
\begin{equation}\label{eq:energy-balance}
    \partial_{t}(\rhoa \ea)+\nabla_{\bx}\cdot(\rhoa \ea\bua +\bqa)+\textnormal{trace}\big(\Pressurea\nabla\bua\big)  = -2\int \limits_{\Omegab(t)} \phiab(\bx,\by)\ea\rho\rhobp\dy, \qquad \bq=\bq_h.
\end{equation}
Here we choose to interpret the equality \eqref{eq:energy-balance} as a special case of  entropic inequality \eqref{eq:meso-pressure}, giving room to validate the  formal passage to the limit in lieu of lack of formal closure.

\medskip\noindent
{\bf The case $p>1$}. We split $(\bv-\bu)\cdot(\bv-\bvp)\equiv  \frac{1}{2}|\bvp-\bv|^2 + \big(\frac{1}{2}(\bvp+\bv)-\bu\big)\cdot(\bvp-\bv) $ to obtain
\begin{equation}\label{eq:split-p}
 \begin{split}
\int \limits_{\R^d} & (\bv-\bu)\cdot \Qp(\fa, \fb)\dv \\
 & = -\int \limits_{{\mathcal S}(t)}\phiab (\bx,\bx')\iint \limits_{\R^d \times \R^d}|\bvp-\bv|^{2p-2}(\bv-\bu)\cdot(\bv-\bvp)\fa \fbp \dvp\dv\dy \\
  & = -\frac{1}{2}\int \limits_{{\mathcal S}(t)}\phiab (\bx,\bx')\iint \limits_{\R^d \times \R^d}|\bvp-\bv|^{2p}\fa \fbp \dvp\dv\dy\\
  & \ \ \ -\int \limits_{{\mathcal S}(t)} \phiab (\bx,\bx')\iint \limits_{\R^d \times \R^d}|\bvp-\bv|^{2p-2}\big(\nicefrac{1}{2}(\bvp+\bv)-\bu\big)\cdot(\bvp-\bv)\fa \fbp \dvp\dv\dy\\
  & :={\mathcal I}_1+{\mathcal I}_2
     \end{split}
  \end{equation}
The the internal integrand in the first term on the right of \eqref{eq:split-p} does not exceed
  \[
  \begin{split}
  -\frac{1}{2}\iint \limits_{\R^d\times \R^d}&|\bvp-\bv|^{2p}\fa \fbp \dvp\dv   \\
   &   \leq   -\frac{1}{2}\Big(\ \iint \limits_{\R^d\times \R^d}|\bvp-\bv|^2\fa \fbp\dvp\dv\Big)^p \Big(\ \iint \limits_{\R^d\times \R^d}\fa \fbp\dvp\dv\Big)^{-\tfrac{p}{p'}} \\
  &\leq  -\frac{1}{2}\Big(\ \iint \limits_{\R^d\times \R^d}|\bvp-\bu'|^2\fa\fbp\dv\dvp + \iint \limits_{\R^d\times \R^d}|\bv-\bu|^2\fa\fbp\dv\dvp\Big)^p(\rho\rhobp)^{-\tfrac{p}{p'}} \\
   &= -\frac{1}{2}(2\rho\rho'\ea+ 2\rho\rhobp\ebp)^p(\rho\rho')^{-\tfrac{p}{p'}} \\
    & \leq  -\hf\big((2\ea)^p+(2\ebp)^p\big)\rho\rhobp.
\end{split}
\]
The first passage on the right follows from  H\"{o}lder inequality, the second   follows from  polarization $\bvp-\bv \equiv (\bvp-\bu')+(\bu'-\bu)+(\bu-\bv)$
and the last from Jensen's inequality. Hence 
\begin{equation}\label{eq:first-term}
{\mathcal I}_1 \leq -\hf \int\limits_{{\mathcal S}(t)} \phiab (\bx,\bx')\big((2\ea)^p+(2\ebp)^p\big)\rho\rhobp.
\end{equation}
For the second term on the right of \eqref{eq:split-p} we claim that it can be written as a complete divergence
  \begin{equation}\label{eq:second-term}
  {\mathcal I}_2(\bx) = \nabla_\bx \cdot \bq_\phi, \qquad 
  \end{equation}
 Indeed, by anti-symmetry $(\bx,\bv)\leftrightarrow (\bxp,\bvp)$  the term ${\mathcal I}_2(\bx)$ has zero mean, 
  \[
  \int \limits_{{\mathcal S}(t)} {\mathcal I}_2(\bx)\dx = \iint \limits_{{\mathcal S}(t)\times {\mathcal S}(t)}\phiab (\bx,\bx')\iint \limits_{\R^d\times \R^d}|\bvp-\bv|^{2p-2}\big(\nicefrac{1}{2}(\bvp+\bv)-\bu\big)\cdot(\bvp-\bv)\fa \fbp \dvp\dv\dx\dy=0.
  \]
 Hence,  there exists a solution , $\Delta \psi = {\mathcal I}_2(\bx), \ \bx\in {\mathcal S}(t)$
  subject to Neumann boundary condition, $\frac{\partial \psi}{\partial {\mathbf n}}_{|\partial{\mathcal S}(t)}=0$, and \eqref{eq:second-term} follows with
   $\bq_\phi=\nabla \psi$.
Combining \eqref{eq:kinternalb} with \eqref{eq:first-term} and \eqref{eq:second-term}  we  arrive at the entrpoic inequality \eqref{eq:meso-ppressure}
\begin{equation}\label{eq:penergy-balance}
\begin{split}
    \partial_{t}(\rhoa \ea)+\nabla_{\bx}\cdot(\rhoa \ea\bua +\bqa)&+\Pressurea\nabla\bua  \\
     & \leq  -\hf\int \limits_{\Omegab(t)} \phiab(\bx,\by)\big((2\ea)^p+(2\ebp)^p\big)\rho\rhobp\dy, \quad \bq:=\bq_h+\bq_\phi.
     \end{split}
\end{equation}
Observe that while the entropic inequality \eqref{eq:energy-balance} in the case $p=1$  was a matter of choice, the corresponding inequality  \eqref{eq:penergy-balance} for  $p>1$ is a matter of necessity in order to make a macroscopic interpretation.

\section{Pointwise bounds on velocity fluctuations}\label{sec:pointwise}

\subsection{Pointwise fluctuations in mono-kinetic alignment\label{sec:pointwise_mono-kinetic}}
Arguing along the lines of  \cite[\S1]{HT2017}, we first fix an arbitrary unit vector ${\bw}\in \mathbb{R}^d$ and project \eqref{eq:hydro-mono-kinetic} onto the space spanned by ${\bw}$ to get
\begin{align*}
(\partial_t+\bua\cdot\nabla_\bx)\lan \bua(t,\bx),{\bw}\ran=\int \limits_{{\mathcal S}(t)}\phi(\bx,\by)(\lan\bu(t,\by),{\bw}\ran-\lan\bu(t,\bx),{\bw}\ran)\rho(t,\by)\dy
\end{align*} 
Now we assume that $\lan\bu(t,\bx),\bw\ran$ reaches  maximum   and  minimum values at $\bx_+=\bx_+(t)$ and, respectively, $\bx_-=\bx_-(t)$,
\[
\begin{split}
u_+(t)&=\lan \bu(t, \bx_+(t)),\bw\ran:=\sup_{\bx\in {\mathcal S}(t)}\lan\bu(t,\bx), \bw\ran \\
 u_-(t)& =\lan \bu_-(t,\bx_+(t)),\bw\ran:=\inf_{\bx\in {\mathcal S}(t)}\lan\bu(t,\bx), \bw\ran.
 \end{split}
\] 
To simplify notations, we temporarily suppress the $\bw$-dependence, $u_\pm(t)=u_\pm(t;\bw)$.
We abbreviate  $\displaystyle \overline{u}(t):=\frac{1}{M}\int \rho\lan \bu(t,\by),\bw\ran\dy$. Since $\lan\bu(t,\by), {\bw}\ran\leq \lan\bu(t,\bx_+),{\bw}\ran$ and, by assumption, $\phi(\bx_+,\by)) \geq k(D(t))$, we find 
\begin{equation}\label{eq:uplus}
\begin{split}
\ddt u_+(t)=&\int \limits_{{\mathcal S}(t)}\phi(\bx_+,\by)\big(\lan\bu(t,\by), {\bw}\ran-\lan\bu(t,\bx_+),{\bw}\ran\big)\rho(t,\by)\dy\\
\leq &k(D(t))\int \limits_{{\mathcal S}(t)} \big(\lan \bu(t,\by),{\bw}\ran-u_+(t)\big)\rho(t,\by)\dy
\\
=&k(D(t))M \big((\ubar(t) -u_+(t)\big).
\end{split}
\end{equation}
 Similarly, we bound
$\displaystyle u_-(t):=\inf_{\bx\in{\mathcal S}}\lan\bu(t,\bx), {\bw}\ran$ obtaining
\[
\ddt u_-(t) \geq k(D(t))M \big(\ubar-u_-(t)\big).
\]
The difference of the last two bounds yields 
\[
\ddt \big(u_+(t)-u_-(t)\big)\leq -k(D(t))M \big(u_+(t)-u_-(t)\big),
\]
and since   $\displaystyle \delbu(t)=\sup_{\bx,\by\in {\mathcal S}(t)}| {\bu}(t,\bx)- {\bu}(t,\by)|=\sup_{|\bw|=1}\big(u_+(t;\bw)-u_-(t;\bw)\big)$ is   the diameter of  velocities projected on arbitrary unit vectors  $\bw$ we end up with
\begin{equation}\label{eq:aprioridelV}
\ddt \delbu(t)\leq -k(D(t))M \delbu(t).
\end{equation}
\subsection{Pointwise fluctuations in mono-kinetic $p$-alignment($p\geq 1$)}\label{sec:pointwise_pmono-kinetic}
We extend the pointwise bound \eqref{eq:aprioridelV} for the general $p$-alignment, $p\geq 1$. By Galilean invariance we may assume $\bm_0=0$, in which case  \eqref{eq:delV-pcontarct} is simplified to  the uniform bound
\begin{equation}\label{eq:max_pmono-kinetic}
\ddt u_+(t) \leq -\frac{1}{2}Mk(D(t))(u_+(t))^{2p-1}, \qquad u_+(t)=\sup_{\bx\in{\mathcal S}(t)}|\bua(\myr{t,\bx)}|.
\end{equation}
Indeed, if we let $\bu_+(t)=\bu(t,\bx_+(t))$ with maximal speed $u_+(t)=|\bu_+(t)|$ along particle path $\ds \myr{\dot{\bx}}_+(t)=\bu(t,\myr{\bx_+(t)})$,   we then find,\footnote{\myr{the precise argument involves  Rademacher lemma, see \cite[Lemma 3.5]{Shv2021}}.} 
\[
\ddt \bu_+(t)= \int \limits_{{\mathcal S}(t)}\phi(\bx,\bxp)|\bu'-\bu|^{2p-2}(\bu'-\bu_+)\rho'\dxp.
\] 
By polarization\myr{,} $\bu_+=\frac{1}{2}(\bu_+-\bubp)+\frac{1}{2}(\bu_++\bubp)$\myr{,} we find
\[
\begin{split}
\frac{1}{2}\ddt |\bu_+(t)|^2 &= -\frac{1}{2}\int \limits_{{\mathcal S}(t)}\phi(\bx,\bxp)|\bu'-\bu_+|^{2p}\rho'\dxp
+ \frac{1}{2}\int \limits_{{\mathcal S}(t)}\phi(\bx,\bxp)|\bu'-\bu_+|^{2p-2}
(|\bu'|^2-|\bu_+|^2)\rho'\dxp \\
& \leq -\frac{1}{2}\int \limits_{{\mathcal S}(t)}\phi(\bx,\bxp)|\bu'-\bu_+|^{2p}\rho'\dxp  \leq -\frac{1}{2}k(D(t))M^{-\tfrac{p}{p'}}\Big(\int \limits_{{\mathcal S}(t)}|\bu'-\bu_+|^{2}\rho'\dxp\Big)^{p} \\
 & \leq -\frac{1}{2}k(D(t))M^pM^{-\tfrac{p}{p'}}|\bu_+(t)|^{2p},
\end{split}
\]
and \eqref{eq:max_pmono-kinetic} follows.
 The first inequality on the right follows from the fact that $|\bu_+|$ is the maximal speed; the second from H\"{o}lder inequality and in the last step we use $\ds \int \bubp\rho'\dxp=0$.
 \subsection{Fluctuations in agent-based description}\label{sec:pagent-based}
 Consider the discrete $p$-alignment model \eqref{eq:palignment} and consider the   energy fluctuations 
\[
\delE(t):=\frac{1}{2N^2}\sum_{i,j=1}^N|\bv_i(t)-\bv_j(t)|^2.
\]
A straightforward computation yields \eqref{eq:pRica}
\[
\begin{split}
\ddt  &\delE(t) = \frac{1}{N^2}\sum_{i,j = 1}^{\Nb}\phiab_{ij}(t)|\bv_{j}(t)-\bv_{i}(t)|^{2p-2}\big\langle \bv_{j}(t)-\bv_{i}(t),\bv_i(t)\big\rangle  \\
 & = -\frac{1}{2N^2}\sum_{i,j = 1}^{\Nb}\phiab_{ij}(t)|\bv_{j}(t)-\bv_{i}(t)|^{2p} 
 \leq -\frac{1}{2}\left(\frac{1}{N^2}\sum_{i,j = 1}^{\Nb}\phiab^{1/p}_{ij}(t)|\bv_{j}(t)-\bv_{i}(t)|^2\right)^p  \\
  & \leq -2^{p-1}k(D(t))\left(\delE(t)\right)^p. 
\end{split}
\]
The first equality follows since $\sum \bv_i(t)$ is conserved in time, and  the second  follows from summation by parts while taking into account the assumed symmetry, $\phi_{ij}=\phi_{ji}$;   next follows the H\"{o}lder inequality (for $p>1$) and finally we use  the lower bound $\phi(\bx_i(t),\bx_j(t))\geq k(D(t))$. 
Similarly, we consider the uniform fluctuations
\[
\delbv(t):=\max_i|\bv_i(t)-\bvbar|.
\]
We assume without loss of generality $\bvbar_0=0 \ \leadsto \bvbar(t)\equiv 0$ and it remains to bound the maximal value $\ds \bv_+(t) = \argmax_{\bv_i}|\bv_i|$. Writing $\bv_+= \frac{1}{2}(\bv_+-\bv_j)+\frac{1}{2}(\bv_++\bv_j)$ we find
\[
\begin{split}
\frac{1}{2}\ddt |\bv_+(t)|^2 &=
\frac{1}{2N}\sum_j \phi_{ij}|\bv_+-\bv_j|^{2p-2} \langle \bv_+-\bv_j, \bv_j-\bv_+\rangle + \frac{1}{2N}\sum_j \phi_{ij}|\bv_+-\bv_j|^{2p-2} \big(|\bv_j|^2-|\bv_+|^2\big)\rangle \\
& \leq  -\frac{1}{2N}\sum_j \phi_{ij}|\bv_+-\bv_j|^{2p} 
\leq -\frac{1}{2N}k(D(t))N^{-\tfrac{p}{p'}}\Big(\sum_j|\bv_+-\bv_j|^2\Big)^p \\
 & \leq -\frac{1}{2N}k(D(t))N^{-\tfrac{p}{p'}}N^p |\bv_+|^{2p},
 \end{split}
\]
and \eqref{eq:pRicb} follows.

\section{From enstrophy bound to H\"{o}lder regularity}\label{sec:GN}
For completeness, we recall here the arguments which lead to the Gagliardo-Nirenberg-Sobolev  inequality, stating that for $\bu\in W^{s,2p}({\mathcal S})$ with `nice' boundary satisfying \ref{H2}, we have H\"{o}lder continuity  of order $s-\frac{d}{2p}$, 
\begin{equation}\label{eq:AGNS}
 |\bu(\bx)-\bu(\by)|\leq C_s\|\bu\|_{{}_{\dot{W}^{s,2p}({\mathcal S})}}|\bx-\by|^{s-\tfrac{d}{2p}}, \qquad \bx,\by\in {\mathcal S}, \ \ \frac{d}{2p} < s<1.
\end{equation}
We follow \cite[theorem 8.2]{DPV2012}. As a first step we note that thanks to hypothesis \ref{H2},  $\bu$ can be extended to $\widetilde{\bu}$ defined over over $\R^d$ with comparable $W^{s,2p}$-norm,  $\|\widetilde{\bu}\|_{W^{s,2p}(\R^d)} \lesssim \|\bu\|_{W^{s,2p}({\mathcal S})}$, 
 \cite[Theorem 5.4]{DPV2012}. We continue with the extension, $\widetilde{\bu}$. Set $R=|\bx-\by|, \ \bx, \by\in {\mathcal S}$ and let $\langle \widetilde{\bu}\rangle_{B_{2R}(\bz)}$ denote the average over the ball $B_{2R}$ centered at $\bz$,
 \[
 \langle \widetilde{\bu}\rangle_{B_{2R}(\bz)}:= \frac{1}{|B_{2R}(\bz)|}\int \limits_{\bz'\in B_{2R}(\bz)} \widetilde{\bu}(\bz'){\d}\bz'.
 \]
 Fix $\bxp'$ as an intermediate point in the intersection of the two balls, $B_{2R}(\bx)\cap B_{2R}(\by)$ and split
 \begin{equation}\label{eq:split-GNS}
 \begin{split}
 |{\bu}(\bx)-{\bu}(\by)|
    \leq |\widetilde{\bu}(\bx) - \langle \widetilde{\bu}\rangle_{B_{2R}(\bx)}|
 &+ |\langle \widetilde{\bu}\rangle_{B_{2R}(\bx)}-\widetilde{\bu}(\bxp')| \\
   &  + |\widetilde{\bu}(\bxp')-\langle \widetilde{\bu}\rangle_{B_{2R}(\by)}|
 + |\langle \widetilde{\bu}\rangle_{B_{2R}(\by)}-\widetilde{\bu}(\by)|. 
  \end{split}
 \end{equation}
By H\"{o}lder inequality, for every $\bw\in B_{2R}(\bz)$, there holds
\begin{equation}\label{eq:local-ave}
\begin{split}
|\langle \widetilde{\bu}& \rangle_{B_{2R}(\bz)}-\widetilde{\bu}(\bw)|  \\
& \leq \frac{1}{|B_{2R}(\bz)|}\int \limits_{\bz'\in B_{2R}(\bz)} |\widetilde{\bu}(\bw)-\widetilde{\bu}(\bz')|{\d}\bz' \leq C_d\|\widetilde{\bu}\|_{W^{s,2p}(B_{2R}(\bz))}R^{s-\tfrac{d}{2p}}, \quad \bw \in B_{2R}(\bz).
\end{split}
\end{equation}
and \eqref{eq:AGNS} follows from proper application of \eqref{eq:local-ave} to each of the terms on the right of \eqref{eq:split-GNS}.\newline 
We close by noting that in the  special 1D case,  \eqref{eq:AGNS}  is reduced to the inequalities of Ladyzheskaya  \cite{MRR2013} or Agmon's \cite[Lemma 13.2]{Agm2010}, 
\[
\max_{\bx\in \Omegaa} |\bu(t,\bx)|\lesssim  \|\bu\|^{1-\tfrac{1}{2s}}_{L^2(\Omega)} \times \|\bu\|_{{H}^s(\Omegaa)}^{\tfrac{1}{2s}}, \quad \nicefrac{1}{2} < s<1.
\]

\section{A uniform dispersion bound}\label{sec:dispersion}
\begin{proof}[Proof\nopunct]  of lemma \ref{lem:dispersion} ({\bf The range $1<p<\nicefrac{3}{2}$}).
We consider the multi-dimensional $p$-alignment, $p>1$, driven by heavy-tailed singular kernel of order $(s,\beta)$. Assume that we have the dispersion bound
\[
\Da(t) \leq C_D(1+t)^{\gamma}, \qquad  \theta=\frac{d}{2p}<s <1.
\]
By \eqref{eq:Dabeta} this holds with $\ds \gamma= \myr{\gamma_p}=\frac{2p-1}{2p(1+\theta-s)}$. We will improve this bound by a bootstrap argument.
To this end we recall that corollary \ref{cor:1} implies the decay (\myr{suffice to consider $t\geq1$})
\[
\delE(t) \leq C_1(\myr{1}+t)^{-\tfrac{1-\beta\gamma}{p-1}}, \qquad \myr{t\geq 1}, \qquad C_1=\myr{2^{\frac{p-1}{1-\beta\gamma}}}C_RM_p.
\]
Here and below, we use the different constants $C_1,C_2, \ldots$ to trace our calculations.\newline
\myr{For the range of $\beta$ assumed in \eqref{eq:Duniform},  
$\beta<(\nicefrac{3}{2}-p)d, \ 1<p<\nicefrac{3}{2}$, we have\footnote{\myr{In fact, the precise bound enables a slightly larger range $\beta < \frac{p(3-2p)}{2p-1}d<(\nicefrac{3}{2}-p)d$ but we prefer to keep it simple with the latter.}}
$\ds 2p-1< \frac{1-\beta\gamma_p}{p-1}$.} 
 Fix $\mu$ such that
 \[
 2p-1<\mu < \frac{1-\beta\gamma_p}{p-1}.
 \]
 The energy fluctuations bound, e.g., \eqref{eq:end-revisited}, yields
\[
\begin{split}
\ddt\Big((1+t)^\mu&\delE(t)\Big) = (1+t)^\mu\ddt \delE(t) + \mu(1+t)^{\mu-1}\delE(t) \\
& \leq -\frac{\rho_-^2}{2}(1+t)^\mu\|\bu(t,\cdot)\|^{2p}_{\dot{W}^{s,2p}({\mathcal S}(t))}
+ C_1\mu(1+t)^{\myr{\mu'}-1}, \qquad \myr{\mu':=\mu-\frac{1-\beta\gamma}{p-1}<0},
\end{split}
\]
and hence  the weighted enstrophy bound
\begin{equation}\label{eq:weighted-ens}
\int \limits_0^t (1+\tau)^\mu\|\bu(\tau,\cdot)\|^{2p}_{W^{s,2p}}{\d}\tau \leq 2C_\rho^2\delE(0)+C_2, \qquad C_2=2C^2_\rho C_1\mu\myr{\frac{1}{|\mu'|}}.
\end{equation}
We now revisit \eqref{eq:dynamicD}, integrating  $\ds \ddt \Da^{1+\theta-s}(t) \leq C'_s\|\bu(t,\cdot)\|_{\dot{W}^{s,2p}({\mathcal S}(t))}$ with a weighted H\"{o}lder inequality, 
 \[
\begin{split}
D^{1+\theta-s}(t) & \leq D_0^{1+\theta-s}
 + C'_s\Big(\int \limits_0^t (1+\tau)^\mu\|\bu(\tau,\cdot)\|^{2p}_{\dot{W}^{s,2p}({\mathcal S}(t))}
{\d}\tau\Big)^{\tfrac{1}{2p}}\Big(\int \limits_0^t (1+\tau)^{-\tfrac{\mu}{2p}(2p)'}\Big)^{\tfrac{1}{(2p)'}}. 
\end{split}
\]
\myr{Using \eqref{eq:weighted-ens} and the fact that $\ds \frac{\mu}{2p-1}>1$ we end up with the uniform bound}
\[
 \myr{D(t) \leq D_+= \big(D_0^{1+\theta-s} + C_3\big)^{\frac{1}{1+\theta-s}}}, \qquad \myr{C_3=C'_s\big(2C_\rho^2\delE(0)+C_2\big)^{\frac{1}{2p}}\Big(\frac{1}{\frac{\mu}{2p-1}-1}\Big)^{\frac{1}{(2p)'}}}
 \]
 \end{proof}
 \begin{remark}[{\bf The case $p=1$}]\label{rem:peq1}
 \myr{It should be possible to extend the uniform dispersion bound of lemma \rem{lem:dispersion} to the limiting case of `pure' alignment, $p=1$. To this end one should use a proper `exponential multiplier', instead of $(1+t)^\mu$ used above for $1<p<\nicefrac{3}{2}$.} 
 \end{remark}
 \begin{remark}[{\bf The case $p>\nicefrac{3}{2}$}]\label{rem:pgt32}
 When $p>\nicefrac{3}{2}$ we are unable to secure a uniform dispersion bound as in  lemma \ref{lem:dispersion}, but we can   still  improve the dispersion rate, $\gamma_p$, using a  more refined bootstrap argument.
For this range of $p$'s we have $\ds \frac{1-\beta\gamma_p}{p-1}< 2p-1$. Fix 
$\mu$ such that $\ds \frac{1-\beta\gamma_p}{p-1}< \mu<2p-1$. 
In this case, $\ds \mu'=\mu-\frac{1-\beta\gamma_p}{p-1}$ is positive and  we have the corresponding enstrophy weighted bound 
\begin{equation}\label{eq:weighted-enspgt32}
\int \limits_0^t (1+\tau)^\mu\|\bu(\tau,\cdot)\|^{2p}_{W^{s,2p}}{\d}\tau \leq 2C_\rho^2\delE(0)+C_2 (1+t)^{\myr{\mu'}}, \qquad C_2=\frac{2C^2_\rho C_1\mu}{\mu'}>0.
\end{equation}
As before, we revisit \eqref{eq:dynamicD},  integrating  $\ds \ddt \Da^{1+\theta-s}(t) \leq C'_s\|\bu(t,\cdot)\|_{\dot{W}^{s,2p}({\mathcal S}(t))}$ with a weighted H\"{o}lder inequality to find
 \[
\begin{split}
D^{1+\theta-s}(t) & \leq D_0^{1+\theta-s}
 + C'_s\Big(\int \limits_0^t (1+\tau)^\mu\|\bu(\tau,\cdot)\|^{2p}_{\dot{W}^{s,2p}({\mathcal S}(t))}
{\d}\tau\Big)^{\tfrac{1}{2p}}\Big(\int \limits_0^t (1+\tau)^{-\tfrac{\mu}{2p}(2p)'}\Big)^{\tfrac{1}{(2p)'}} \\ 
 & \leq D_0^{1+\theta-s} + C_3(1+t)^{\frac{\mu'}{2p}}\times (1+t)^{\big(1-\frac{\mu}{2p-1}\big)\frac{1}{(2p)'}} \\
 & = D_0^{1+\theta-s} + C_3(1+t)^{\tfrac{1}{2p}\Big((2p-1)-\frac{1-\beta\gamma}{p-1}\Big)},
 \qquad C_3=C'_sC_2\Big(1-\frac{\mu}{2p-1}\Big)^{-\tfrac{1}{(2p)'}}>0.
\end{split}
\]
We conclude with a dispersion bound
\begin{equation}\label{eq:gamma-iterations}
D(t) \leq C'_D(1+t)^{\gamma'}, \qquad \gamma':=\frac{1}{2p(1+\theta-s)}\Big((2p-1)-\frac{1-\beta\gamma}{p-1}\Big)< \gamma_p.
\end{equation}
Recall that the requirement $\beta\gamma_p<1$ led to the $\beta$-restriction
in \eqref{eq:mDheavys}, $\ds \beta<\frac{d}{2p-1}$. Therefore, 
\[
\frac{\beta}{2p(1+\theta-s)(p-1)} <\frac{\beta}{d(p-1)} <\frac{1}{(2p-1)(p-1)}<1, \qquad  p> \nicefrac{3}{2},
\]
 so that the fixed point iterations in \eqref{eq:gamma-iterations}, $\gamma \mapsto \gamma'$, contract towards a limiting value $\gamma=\gamma_*$
\[
\gamma_*=\frac{p(2p-3)}{(p-1)2p(1+\theta-s)-\beta},\qquad p> \nicefrac{3}{2}, \quad \beta< \frac{d}{2p-1}.
\]
In particular, since $2p(1+\theta-s)>d$, then after finitely many iterations there holds 
\begin{equation}\label{eq:gammagt32}
D(t) \leq C'_D(1+t)^\gamma, \quad \gamma= \frac{2p\big(p-\nicefrac{3}{2}\big)}{(p-1)d-\beta},   \qquad \beta < \frac{d}{2p-1}, \quad  p>\nicefrac{3}{2},
\end{equation}
which improves the bound of $\ds \frac{2p-1}{d} > \gamma_p$.
\end{remark}


\end{document}